\documentclass[11pt,a4paper,psamsfonts,twoside]{amsart}
\usepackage{latexsym, amssymb,amscd,amsmath,amsthm}
\usepackage{graphicx,color} 

\chardef\bslash=`\\ 

\newtheorem{theorem}{Theorem}[section]
\newtheorem{corollary}[theorem]{Corollary}
\newtheorem{lemma}[theorem]{Lemma}
\newtheorem{proposition}[theorem]{Proposition}

\theoremstyle{remark}
\newtheorem{remark}[theorem]{Remark}
\newtheorem{example}[theorem]{Example}
\theoremstyle{definition}
\newtheorem{definition}[theorem]{Definition}
\newtheorem*{notation}{Notation}

\numberwithin{equation}{section}

\newcommand{\thmref}[1]{Theorem~\ref{#1}}
\newcommand{\secref}[1]{Section \ref{#1}}
\newcommand{\proref}[1]{Proposition~\ref{#1}}
\newcommand{\lemref}[1]{Lemma~\ref{#1}}
\newcommand{\corref}[1]{Corollary~\ref{#1}}

\def\Cal#1{{\mathcal#1}}
\def\<{\langle}\def\>{\rangle}
\def\leq{\leqslant}
\def\geq{\geqslant}
\def\what{\widehat}
\def\Z{{\mathbb Z}}\def\N{{\mathbb N}} 
 
\def\X{{\bf X}}

\def\Spec{\textsl{Spec}}
\def\Fix{\textsl{Fix}}

\def\restr{\upharpoonright}
\newcommand{\inv}{^{-1}}

\def\hsthree{\quad}
\def\hsfour{\quad}

\def\hsseven{\qquad}
\def\hsten{\qquad\quad}
\def\bsl{\backslash}

\def\fin{  }

\def\al{\alpha}                 \def\be{\beta}		\def\lam{\lambda}
                 \def\Ga{\Gamma}		\def\Lam{\Lambda}

               \def\varep{\varepsilon}
\begin{document}

\title[Boundary quotients and ideals of Toeplitz C*-algebras]
{Boundary quotients and ideals of \\
Toeplitz C*-algebras of Artin Groups}
\date{1 April  2006}
\author[J.~Crisp]{John Crisp$^1$}
	\address{I.M.B.(UMR 5584 du CNRS)\\ 
              Universit\'e de Bourgogne\\
		     B.P. 47 870\\
		21078 Dijon, France}

\author[M.~Laca]{Marcelo Laca$^2$}
\address{Department of Mathematics and Statistics
       University of Victoria, Victoria, BC, V8W 3P4, CANADA}

\thanks{$^1)$ Supported by the CNRS, France.}
\thanks{$^2)$ Supported by the NSERC, Canada.}

\keywords{quasi-lattice order, covariant
isometric representation, Toeplitz algebra, Artin group.}


\begin{abstract} 
We study the quotients of the Toeplitz C*-algebra of a quasi-lattice ordered group $(G,P)$,
 which we view as crossed products by a partial actions of $G$ on closed invariant subsets of 
a totally disconnected compact Hausdorff space, the Nica spectrum of $(G,P)$.
Our original motivation and our main examples
 are drawn from right-angled Artin groups,
 but many of our results are valid for more general quasi-lattice ordered groups.
We show that the Nica spectrum has a unique minimal closed invariant subset, 
which we call the boundary spectrum, and we define the boundary quotient to be
the crossed product of the corresponding restricted partial action.
The main technical tools used are the results of Exel, Laca, and Quigg on simplicity and 
ideal structure of partial crossed products, which depend on amenability 
and topological freeness of the partial action and its restriction to closed invariant subsets.
When there exists a generalised length function, or controlled map, defined
on $G$ and taking values in an amenable group, we prove that
 the partial action is amenable on arbitrary closed invariant subsets. 
 The topological freeness of the boundary action depends on topological freeness
of the restriction to a certain lattice subgroup of $G$,  the ``core" of
$(G,P)$, which often turns out to be trivial.
Our main results are obtained for right-angled Artin groups with trivial centre, 
that is, those with no cyclic direct factor; they include
a presentation of the boundary quotient in terms of generators and relations that 
generalises Cuntz's presentation of  $\mathcal O_n$, 
a proof that the boundary quotient  is purely infinite and simple, and a parametrisation of the ideals of the 
Toeplitz C*-algebra in terms of subsets of the standard generators of the Artin group. 
\end{abstract}

\maketitle

\section*{Introduction}

To any quasi-lattice ordered group $(G,P)$, one can associate the Toeplitz  (or Wiener-Hopf)
C*-algebra $\Cal T(G,P)$ generated by the compression to $\ell^2(P)$  of the left regular representation of $G$. 
These were originally introduced by Nica  \cite{nica}, and were further studied in  \cite{quasilat,ELQ}, in terms of 
a characteristic universal property of representations of the semigroup $P$ by isometries satisfying a certain 
covariance relation isolated by Nica. Along these lines, it was shown in \cite{CL} that
the Toeplitz algebras associated in this way to 
right angled Artin groups have the universal property with
respect to covariant isometric representations and satisfy 
a uniqueness theorem, \cite[Theorem 24]{CL}, that
generalises results of Coburn \cite{cob} and Cuntz \cite{cun2}.
A key feature of this uniqueness result is that only 
the isometric representations of the Artin monoid that 
satisfy a certain properness condition, see  \eqref{propercond} below,
 give rise to faithful representations of the Toeplitz algebra.
 In the case of Coburn's classical result, concerning a single isometry, 
this properness condition simply says that the isometry must not be a unitary, while in 
Cuntz's situation, where the algebra is generated by a number 
of isometries with mutually orthogonal ranges, it says that the direct sum of these ranges is
not the whole space.
 
It is then natural to consider the different ways in which the properness condition may fail to hold, 
and the ideals and quotients of the Toeplitz C*-algebras that arise from this failure.
We continue here the study of these questions, the importance of which is underlined by the
well known facts that
 in Coburn's situation one obtains the C*-algebra of continuous functions on the circle
(and its quotients), and in Cuntz's situation one obtains the purely infinite 
Cuntz algebras $\mathcal O_n$ from \cite{cun}. Our present motivation is to understand the situation for 
right-angled Artin groups, the study of which was initiated in \cite{CL}, 
but as it turns out, many of our results are valid  for more general classes of
quasi-lattice ordered groups, notably those which admit a certain type of 
generalized length function, or ``controlled map'', taking values in an amenable group.

The methods used for the analysis are those developed in \cite{ELQ} and involve realising
 Toeplitz algebras as crossed products by partial actions.
Once topological freeness and amenability have been established for
the relevant partial actions, these methods  allow us to move between ideals of
the Toeplitz C*-algebras and closed invariant subsets of the spectrum of the diagonal,
which, as shown by Nica in \cite{nica}, is the space of
nonempty, hereditary, directed subsets of  the positive cone $P$ of the quasi-lattice ordered
group $(G,P)$ under consideration.

We begin by giving in \secref{Sect1} a
 brief reminder of  the relevant results about partial actions from \cite{ELQ} stated 
with a small modification that will be needed here. In \secref{Sect2} 
we  include some basic facts about quasi-lattice ordered groups and 
their Toeplitz C*-algebras, and 
also a brief review of the definition of right-angled Artin groups.

In \secref{SectNica} we discuss closed invariant subsets of the Nica spectrum and the
induced ideals of the partial crossed product they generate.
Parallel to this, we also explore the different
  ways in which the properness condition \eqref{propercond} 
can fail for covariant isometric representations. 
Of particular interest are two collections of extra relations 
that can be imposed on the isometric representations. 
One is a maximal set of relations $\Cal F$ which corresponds to the closure 
of the set of maximal points in the Nica spectrum, and determines 
the boundary quotient. The other set $\Cal E$ is empty unless 
$P$ admits a finite set of lower bounds, in which case it corresponds 
to the set of unbounded points of the Nica spectrum 
and determines the quotient by a minimal ideal.
  
  Amenability is discussed in \secref{SectAmenable}. 
 The main results there are \thmref{amenable} and  \corref{coramenable},
which apply to quasi-lattice ordered groups admitting a controlled map to 
an amenable group. By extending the methods of \cite{quasilat} we are able 
to give a direct proof of  amenability of the partial action when
restricted to arbitrary closed invariant subsets
of the Nica spectrum. This argument avoids using the approximation property 
of \cite{ELQ}, which we have been unable to establish in the case of Artin groups.

 In \secref{SectBoundary}  we focus on the 
partial action restricted to the  boundary. \thmref{PureSimple} states that the boundary
quotient is purely infinite and simple whenever this boundary action is both amenable
and topologically free.
In order to characterise topological freeness in \proref{TopFreeCore} 
we introduce the notion of the \emph{core} of a quasi-lattice ordered group. 
For Artin groups this corresponds to the centre. When it is trivial, 
the partial action of the right-angled Artin group 
on its boundary spectrum  is topologically free, by \corref{TopFreeArtin}.
 
In \secref{SectIrredCase} we introduce the notion of graph-irreducibility 
for a quasi-lattice ordered group $(G,P)$ and show, \proref{IrredCase},
that if $(G,P)$ admits a controlled map to a free abelian group and is graph-irreducible
then the boundary relations $\Cal F$ and minimal set of relations $\Cal E$ are 
equivalent, and so the boundary ideal is either minimal or trivial. 
In the case of a right-angled Artin group, graph-irreducibility coincides with the 
typically weaker notion of irreducibility with respect 
to direct sums (of partially ordered groups). Moreover the decomposition
of a right-angled Artin group $(A,A^+)$ into a direct sum of (graph-)irreducibles 
can be easily read from the usual presentation for the group -- the presentation 
is described by a graph $\Gamma$ and the irreducible factors correspond to the
connected components of the opposite graph $\Gamma^{\rm opp}$ 
(by definition, $\Gamma^{\rm opp}$ has the same vertex set as $\Gamma$
and edges joining the vertices that are not joined in $\Gamma$).

These considerations lead us to our first main result, 
 \thmref{presentationboundaryquotient}, in which we prove that, 
  for each right-angled Artin group with trivial centre, the boundary quotient
of the Toeplitz algebra is purely infinite and simple and has 
 a straightforward presentation in terms of generators and relations.
This generalizes Cuntz's classical result for $\mathcal O_n$, which is
the boundary quotient associated to the free group 
on $n$ generators. Given the richness of the class of right-angled Artin groups,
this result raises the interesting question of classification of thes boundary quotients, 
which is not addressed in the present work, but is the subject of ongoing 
joint research with B. Abadie.
 
 In \secref{SectProducts} we study the decomposition of $(G,P)$ as a 
direct sum of quasi-lattice ordered groups.  The motivating example is the
 direct sum decomposition of a right-angled Artin group according to the 
connected components of $\Gamma^{\rm opp}$.
We show that such a sum is topologically free on closed 
invariant subsets if and only if all the summands are. Finally, 
in \secref{SectIdealStruct}  we combine this together with the results of 
\secref{SectAmenable}, \secref{SectBoundary}, and \secref{SectIrredCase}, and 
the parametrisation of ideals given in \cite{ELQ} to obtain our second main result,  
\thmref{LatticeOfIdeals},  which shows that the ideals of the Toeplitz C*-algebra
of a right-angled Artin group are parametrised by the elements of the Boolean algebra of 
finite subsets of the set of all finite connected components of $\Gamma^{\rm opp}$.


\section{Partial group actions and crossed products}\label{Sect1} 


Let $\X$ denote a locally compact topological space. By a \emph{partial action}
$\theta$ of a group $G$ on $\X$ we mean a family of open sets $\{ U_t : t\in G\}$
and partial homeomorphisms $\theta_t:U_{t^{-1}}\to U_t$ such that $\theta_{st}$
extends $\theta_s\theta_t$ for all $s,t\in G$ (see \cite{McCl,Exel2,ELQ}). Such an action 
naturally induces a partial action $\al$ of $G$ on the C*-algebra $C_0(\X)$ with domains 
the ideals $D_t:=C_0(U_t)$ and partial isomorphisms 
$\alpha_t : f\in D_t \mapsto f\circ \theta_t\in D_{t\inv}$
for $t\in G$. We refer to the triple $(C_0(\X), G, \al)$
as a \emph{partial dynamical system}. 
Associated to any partial dynamical system
$(C_0(\X), G, \al)$ there is a crossed product $C^*$-algebra
$
C_0(\X)\rtimes_\al G
$
and a reduced  crossed product $C^*$-algebra
$C_0(\X)\rtimes_{\al,r} G$. By analogy with the corresponding 
constructions for group actions, these are defined as  completions of the convolution 
algebra of finite sums of the form $\sum f_g \delta_g$ where $f_g\in D_g$ for each $g$.
The full crossed product $C_0(\X)\rtimes_\al G$ can also be defined
 in terms of a universal property for covariant representations.
See \cite{McCl,qui-rae,Exel2,ELQ} for details. \fin

In order to study simplicity, pure infiniteness and the 
ideal structure of crossed products by partial actions 
we rely heavily on the methods developed in \cite{ELQ}. 
 These methods depend upon
the following two key properties which we shall therefore need to verify 
for each of our partial dynamical systems. 

\begin{definition} 
We say that a partial action of a group $G$ on a locally compact 
topological space $\X$, or the induced partial action on the 
C*-algebra $C_0(\X)$, is \emph{topologically free} if the fixed 
set in $\X$ of every nontrivial element of $G$ has empty interior.

We say that a partial action is \emph{amenable} if the canonical 
map from the full to the reduced crossed product is an isomorphism, 
equivalently, if the conditional expectation 
$\Phi : C_0(\X)\rtimes_\al G\to C_0(\X)$,
characterised by 
$\Phi(f\delta_g) = f$ if $g=e$ and $0$ otherwise,
is faithful on positive elements.
\end{definition}

The main result needed about crossed products by partial actions 
is \cite[Theorem 3.5]{ELQ}.  We point out that the approximation property 
assumed there is only needed to prove amenability on closed invariant subsets.
Since we shall verify directly that our partial actions are amenable on closed 
invariant subsets, the following restatement of 
\cite[Theorem 3.5]{ELQ} will be more useful.

\begin{theorem}\label{IdealStructure}
Let  $(C_0(\X), G, \al)$ be a partial dynamical system that is topologically
free and amenable on every closed invariant subset of $\X$. Denote by $\< S\>$
the ideal generated by a set $S$. Then the map
\[
U\mapsto \< C_0(U)\>
\]
is an isomorphism between the lattice of invariant open subsets of $\X$ and the lattice
of ideals of the crossed product $C_0(\X)\rtimes_\al G$.
Moreover, the quotient by the ideal $\<C_0(U)\>$ is canonically
isomorphic to $C_0(\X\setminus U)\rtimes_\al G$.
\end{theorem}

The class of partial actions we shall consider arise as in \cite{ELQ}.
If $G$ is a countable group, with identity element $e\in G$, we define the space
\[
\X_G=\{\omega\in\{0,1\}^G : e\in\omega\}\,.
\]
This is a compact Hausdorff space with the relative topology 
inherited from $\{0,1\}^G $. There is a canonical partial action of $G$ on $\X_G$
(by partial homeomorphisms) which is defined by left multiplication: for $t\in G$, set 
$U_t=\{ \omega\in \X_G:t\in\omega\}$ and define the partial action by
\[
t:U_{t^{-1}}\to U_t\hsthree\text{ such that } \quad 
\omega\mapsto t\omega=\{ tx: x\in\omega\}\,.
\] 

Covariant representations of this partial C*-dynamical system are in 
one to one correspondence with partial representations of the group $G$, 
see \cite{Exel2,ELQ}. Thus the crossed product $C_0(\X)\rtimes_\al G$ 
has a universal property with respect to partial representations 
of $G$.  As in \cite{ELQ}, we are interested here in
partial representations subject to relations, so we adopt the same definitions.
We shall regard a collection $\Cal R$ of continuous functions on $\X_G$ as \emph{relations}
and define the \emph{spectrum $\Spec(\Cal R)$} of $\Cal R$ to be the subset of $\X_G$
\[
\Spec(\Cal R) := \{ \omega \in \X_G : f(t^{-1}\omega)=0 
\text{ for all } t\in \omega, f\in \Cal R\}
\] 
Note that $\Spec(\Cal R)$ is the largest  subset of $\X_G$ on which 
the equivariant form of the relations are satisfied; it is a closed, hence compact set,
and is invariant under the partial 
action of $G$ even if the given set of relations is not $G$-invariant to begin with.
Thus there is a restriction of the partial action to $\Spec(\Cal R)$ 
 and we may form the crossed product $C(\Spec(\Cal R))\rtimes G$. 
On the other hand, every closed invariant subset $X$ of $\X_G$ 
may be obtained as the spectrum of a set of relations 
(namely, those functions in $C(\X_G)$ which are zero on $X$). 
It is shown in \cite[Proposition 4.1]{ELQ} that if $I=\<\Cal R\>$ is the ideal of 
$C(\X_G)$ generated by the set of relations $\Cal R$ then 
$I=C_0(\X_G\setminus\Spec(\Cal R))$ and the quotient is canonically isomorphic to
$C(\Spec(\Cal R))$. Moreover, one also has the following 
exact sequence from \cite[Proposition 3.1]{ELQ}),
\[
0\to \<\Cal R\>\rtimes G \to C_0(\X_G)\rtimes G\to C(\Spec(\Cal R))\rtimes G \to 0 \,,
\]
where the latter crossed product is universal for partial representations of $G$
subject to the relations $\Cal R$, in the sense of \cite[Theorem 4.4]{ELQ}. 
Since $\Spec(\Cal R)$ is chosen so that the
relations in $\Cal R$ are satisfied pointwise on each $\omega \in \Spec(\Cal R)$,  we will
often reverse the terminology, and say that ``an element $\omega$ of $X_G$ 
satisfies the relations $\Cal R$" whenever $\omega \in \Spec(\Cal R)$. 
Obviously the partial representation of $G$ arising in the crossed product 
$C(\Spec(\Cal R))\rtimes G$
satisfies the relations $\Cal R$, so these relations are
imposed via restriction of the partial action from $\X_G$ to $\Spec(\Cal R)$.

\section{Quasi-lattice orders and their Toeplitz algebras}\label{Sect2}

We are interested in certain crossed products arising from quasi-lattice 
ordered groups, as introduced by  Nica in \cite{nica}, see 
also \cite{quasilat,ELQ, CL}. We briefly review the basic facts about these structures.

Let $G$ be a group, with identity element $e$, and suppose that $G$ is equipped
with a partial order $\leq$ that is invariant by left multiplication in the group: 
$x\leq y\implies gx\leq gy$ for all $g,x,y\in G$.
The \emph{positive cone} of $(G,\leq)$ is defined to be the set $P=\{g\in G : e\leq g\}$.
By left invariance of the partial order one has that
$x\leq y$ if and only if $x^{-1}y\in P$, for $x,y\in G$. Thus any left 
invariant partial order on a group is uniquely determine 
by its positive cone.
Moreover, one easily checks that 
a subset $P\subset G$ is the positive cone of a left invariant partial order
on $G$ if and only if  $P$ is a submonoid of $G$ and $P\cap P^{-1}=\{ e\}$. 
In this case we refer to the pair $(G,P)$ as a \emph{partially ordered group}
and denote by $\leq$ the associated left invariant partial order. 
(Note that there is a similarly defined unique
\emph{right} invariant partial order canonically associated to $(G,P)$. However,
in the present paper, we shall consider only left invariant structures).

\begin{definition}
A partially ordered group $(G,P)$ is said to be \emph{quasi-lattice ordered}
if every pair of elements having a common upper bound in $G$ has a least common upper bound 
(with respect to the left partial order $\leq$). The least upper bound
of the pair $x,y\in G$ is usually denoted $x\vee y$. We note that in a
quasi-lattice order every finite set $F$ of elements 
with a common upper bound has a least upper bound, written $\bigvee F$. We
shall generally also write $x\vee y=\infty$ to mean that $x$ and $y$ have no common 
upper bound in $(G,P)$.
\end{definition}

A semigroup  representation $V: P \to \mathcal B(\mathcal H)$ of $P$ by isometries on the Hilbert 
space $\mathcal H$ satisfying $V_xV_x^* V_yV_y^* = V_{x\vee y}V_{x\vee y}^*$ for
$x, y\in P$,  is said to be covariant in the sense of Nica,  \cite{nica}. 
The universal C*-algebra $C^*(G,P)$ for  Nica-covariant isometric representations
of $P$ is canonically isomorphic to the crossed product associated 
to a partial action of $G$ as described in Section \ref{Sect1},
while the Toeplitz C*-algebra
$\mathcal T(G,P)$  determined by the left regular representation of
$P$ on $\ell^2(P)$ is canonically isomorphic to the corresponding reduced crossed product,
see \cite[Section 6]{ELQ}. 
The covariance relations above
 translate into relations for the canonical partial action of $G$ on $\X_G$ given by
\[
 \Cal N := \{ 1_x1_y - 1_{x\vee y} :  \text{ for all } x,y\in G\},
\]
where $1_x$ denotes the characteristic function of the set $\{\omega \in X_G: x\in \omega\}$.
Notice that the extra relation $u_x^* u_x = 1$ for $x\in P$ listed in \cite{ELQ} 
is in fact a consequence of the above because for every $x\in P$ one has that 
$x\inv \vee e = e$ and thus $u_x^* u_x = 1_{x\inv}= 1_{x\inv}1_e= 1_{x\inv \vee e} = 1_e = 1$.

A subset $\omega$ of $G$ is {\em hereditary} if 
$\omega P^{-1}\subset\omega$, and it is {\em directed} if $x\vee y\in\omega$ 
(in particular $x\vee y\neq\infty$) for all $x,y\in\omega$. 
From \cite[\S 6.2]{nica}, the spectrum of the relations $\Cal N$ is the set
\[
\Omega = \Spec(\Cal N) = \{ \text{ directed hereditary subsets of $G$ that contain $e$}\}\,.
\]
As remarked above, we have canonical isomorphisms
\[
C^*_u(G,P)\cong C(\Omega)\rtimes G\quad \text{ and }\quad 
\Cal T(G,P)\cong C(\Omega)\rtimes_r G\,.
\]

When these two algebras are canonically isomorphic,
that is, when the partial action on $\Omega$ is amenable, \thmref{IdealStructure}
describes the  lattice of ideals, but in
 order to use it profitably, we need to determine first  the 
closed invariant subsets of $\Omega$, and then to
check whether the restricted partial actions are amenable and topologically 
free on each of these sets.

These issues will be the focus for the remainder 
of the paper. Some general observations concerning 
closed invariant sets are made in \secref{SectNica}, below. 
The problem of determining all
closed invariant sets is addressed more fully in \secref{SectIrredCase} and \secref{SectProducts}, 
and is concluded in the case of a right-angled Artin group in \secref{SectIdealStruct}.
The question of amenability is addressed in \secref{SectAmenable}, and topological freeness
is studied in \secref{SectBoundary}, and then further in \secref{SectProducts}. 
 
 We conclude this section by describing our main example of a quasi-lattice ordered group. 

\begin{definition}[Right-angled Artin groups]
Let $\Gamma$ denote a simplicial graph with countable vertex set $S$. We define
the corresponding \emph{right-angled Artin group} $A_\Gamma$ by the presentation
\[
A_\Gamma =\<\, S\ \mid\ ab=ba \text{ if } \{ a,b\} \text{ is an edge of }\Gamma\,\>\,.
\]
Let $A_\Gamma^+$ denote the submonoid of $A_\Gamma$ generated by $S$. Then
$(A_\Gamma,A_\Gamma^+)$ is a quasi-lattice ordered group where, for distinct generators
$a,b\in S$ we have $a\vee b=ab=ba$ if $a$ and $b$ span an edge and $a\vee b=\infty$
otherwise (see \cite{CL}).  It is important to note also that $A_\Gamma^+$ is 
presented abstractly as a monoid by the same presentation over $S$ as just 
given for the group $A_\Gamma$. 
This fact is an easy consequence of the 
solution to the word problem in a right-angled Artin group
which seems to have first appeared in the paper \cite{Baudisch} by Baudisch. 
See also \cite{CW}, Section 5, for a more recent treatment of the word and 
conjugacy problems, and for further references to the literature. 
\end{definition}

The Toeplitz C*-algebras of right-angled Artin groups were studied in \cite{CL}, 
where it was shown that a covariant isometric representation $V$ of a right-angled 
Artin monoid gives a faithful representation of the associated Toeplitz C*-algebra 
if and only if  
\begin{equation}\label{propercond}
\prod_{s\in F} (I - V_s V_s^*) \neq 0 \qquad 
\text{ for every finite subset $F$ of generators. }
\end{equation} 
When the graph defining the Artin group consists of vertices with no edges, 
this retrieves Cuntz's result on the uniqueness of the C*-algebra generated 
by isometries with orthogonal ranges that do not add up to the whole space. 
At the other extreme, the case of a full graph yields a result on the uniqueness 
of the C*-algebra generated by n-tuples of $*$-commuting isometries \cite{cob,sal}.

\section{Closed invariant sets in the Nica spectrum}\label{SectNica}
Throughout this section we suppose that $(G,P)$ is a quasi-lattice ordered group.
We wish to study the closed invariant subsets of the Nica spectrum 
$\Omega= \Spec(\Cal N)$  associated to $(G,P)$. 
By the discussion in Section \ref{Sect1}, every closed invariant 
subset of $\Omega$ can be obtained by adding further relations 
to $\Cal N$. From now on, whenever $(G,P)$ is quasi-lattice ordered, 
the Nica relations are implicitly assumed, and only the 
extra relations are indicated explicitly.  We shall see below 
that since the topology on $\Omega$ arises from the order 
structure of $(G,P)$, it will suffice to consider relations 
of a very simple kind, coming from finite subsets of $P$.

\begin{definition} 
Let $\Cal Q$ denote the collection of all finite subsets of $P$. 
The {\em elementary relations} are those of the form 
$f_H := \prod_{h\in H} (1 - 1_h)$ with $H\in\Cal Q$. We
include here the empty set $\emptyset\in\Cal Q$, the 
empty product $f_\emptyset$ being interpreted as $1$.
\end{definition}

For simplicity of notation, we shall frequently refer to the finite 
subsets of $P$, that is, to elements of $\Cal Q$ themselves, as elementary 
relations, and shall say that $\omega\in\Omega$ 
\emph{satisfies the elementary relation $H\in\Cal Q$} 
if the relation $f_H$ is satisfied at $\omega$, i.e. if 
$f_H (t\inv \omega) = \prod_{h\in H} (1 - 1_h)(t\inv \omega) = 0$ 
for every $t\in \omega$.
Using the definition of $1_h$ we see that this means 
that  $\omega\cap tH\neq\emptyset$ for all $t\in\omega$.
Given a set of elementary relations $\Cal R\subset\Cal Q$, 
we shall write $\Omega_{\Cal R}=\Spec(\Cal R)\cap\Omega$ for the 
spectrum of the Nica relations together with the relations in $\Cal R$. 
Thus
\[
\begin{aligned}
\Omega_{\Cal R} &= \{\ \omega \in \Omega\ :\  f_H(t\inv \omega) = 0 
\text{ for all } t \in \omega \text{ and } H\in \Cal R\ \}\\
&= \{\ \omega\in \Omega\ :\ \omega\cap tH\neq\emptyset
\text{ for all } t \in \omega \text{ and } H\in \Cal R\ \}\,.\\
\end{aligned}
\]
We remark that the elementary relations $\Cal Q$ encompass two extreme examples.
The empty set $\emptyset\in\Cal Q$, interpreted as the relation $f_\emptyset =1$, 
is satisfied by \emph{no element} of $\Omega$.
 On the other hand, an elementary relation 
$H\in\Cal Q$ is satisfied by \emph{every} element of $\Omega$ if and only if $e\in H$.

By definition, the topology on $\Omega$ is inherited from the product topology on $\{0,1\}^G$, 
but it also has another characterisation in terms of the 
order structure on $P$, cf. \cite{purelinf}. According to this characterisation, 
a basis of the topology is given by the family of clopen subsets of $\Omega$
\[
V(t,tH) := \{\ \omega\in\Omega\ :\ t\in\omega \text{ but } \omega\cap tH =\emptyset\ \}
\quad \text{ for } t\in G \text{ and } H\in\Cal Q. 
\]
Note that the basic sets $V(t,t\emptyset)$ are included, since $\emptyset\in\Cal Q$. 
These are necessary: in the case where $(G,P)$
is a lattice, the group $G$ itself is a directed hereditary set, so 
$G\in \Omega$, and it is easy to check that  $G\in V(t,tH)$
if and only if $H=\emptyset$.
 
The following lemma shows that in the study of closed invariant subsets of 
$\Omega$, there is no loss of generality in considering only
elementary relations.

\begin{lemma}\label{ClosedSetsRelations}
Every closed invariant subset of the Nica 
spectrum may be written as $\Omega_{\Cal R}$ for
some set of elementary relations $\Cal R\subset\Cal Q$.
More precisely, given any subset $X\subset \Omega$, the smallest closed 
invariant subset of $\Omega$ containing $X$ is equal to $\Omega_{\Cal L(X)}$
where $\Cal L(X)$ denotes the largest collection of 
elementary relations that are satisfied everywhere on $X$, namely
\[
\begin{aligned}
\Cal L(X)&=\{ H\in\Cal Q\ :\ \omega\cap tH \neq\emptyset 
\text{ for all }\, t\in\omega \text{ and all }\, \omega\in X\ \}\\
&=\{ H\in\Cal Q\ :\  X\subset \Omega_{\{ H\}}\,\}\,.
\end{aligned}
\]
\end{lemma}

\begin{proof}
 For $H\in\Cal Q$, let $V(H)=\Omega\setminus \Omega_{\{H\}}=
\{\, \omega\in\Omega\, :\, \omega\cap tH=\emptyset \text{ for some } t\in\omega\, \}$, 
and observe that $V(H)=\bigcup\limits_{t\in G} V(t,tH)$. 
Since the open sets $V(t,tH)$ form a basis for the topology, any open invariant
set can be covered by sets $V(H)$ for $H\in\Cal Q$. In other words, any closed 
invariant set has the form 
$\Omega\setminus \bigcup\limits_{H\in\Cal R} V(H)=
\bigcap\limits_{H\in\Cal R} \Omega_{\{ H\}} =\Omega_\Cal R$
for some $\Cal R\subset\Cal Q$. 
Finally, given any set $X\subset \Omega$, it is clear that 
$X\subset\Omega_{\Cal R}$ if and only if $\Cal R\subset\Cal L(X)$.
Thus $\Omega_{\Cal L(X)}$ is the smallest closed invariant set containing $X$.
\end{proof}
Of course,  when $X$ is already closed and invariant then $\Omega_{\Cal L(X)}=X$.

\begin{definition}
We shall say that a set  of elementary 
relations $\Cal R\subset\Cal Q$ is \emph{saturated} if $\Cal L(\Omega_{\Cal R})=\Cal R$.
\end{definition}

\begin{definition}\label{boundarydef}
We define the \emph{boundary of the Nica spectrum}, or simply the 
\emph{boundary spectrum},  to be the 
spectrum $\Omega_\Cal F$ of the collection of elementary relations
\[
\Cal F := \{\ F\in\Cal Q\ :\ \text{for all } z\in P \text{ there exists } 
x\in F \text{ such that } x\vee z\neq\infty\ \}\,.
\]
Note that $\emptyset\notin\Cal F$. 
\end{definition}

A different definition of the boundary $\partial \Omega$ was given 
 in \cite{purelinf}, in terms of maximal directed hereditary subsets of $G$,
and we need to show that the two definitions are equivalent, that is, $\Omega_{\mathcal F} = \partial \Omega$. 
The following lemma is 
essentially from \cite{purelinf}; the relations $\Cal F$ are not mentioned explicitly there
but they  appear implicitly in \cite[Lemma 5.2]{purelinf}.
We note in passing that 
a maximal directed subset of $G$ is automatically hereditary, because
whenever $\omega\subset G$ is directed, its hereditary closure $\omega P\inv$ is 
also directed and contains $\omega$.

\begin{lemma} \label{minimalityofboundary}
The boundary spectrum $\Omega_\Cal F$ is the (unique) smallest nontrivial closed invariant 
subset of the Nica spectrum $\Omega$, and is equal to the closure $\partial \Omega$ of the 
set of maximal directed subsets of $G$ containing the identity.
\end{lemma}

\begin{proof}
That $\partial\Omega$ is the unique minimal nontrivial 
closed invariant subset of $\Omega$ is proved in \cite[Theorem 3.7]{purelinf}. 
We shall show the equality $\Omega_{\Cal F}=\partial\Omega$. 

Suppose that $K$ is a finite subset of $P$ that is not in $\Cal F$. 
Then there exists $y\in P$ such that
$x\vee y=\infty$, for all $x\in K$. By directedness this implies that
$y\omega\cap K=\emptyset$ for all $\omega\in\Omega$. Thus
every $\omega\in\Omega$ fails to satisfy the relation $K$ (at the value $y^{-1}\in\omega$).
In other words, any set of elementary relations containing $K$ has empty spectrum (and 
therefore has saturation $\Cal L(\emptyset)=\Cal Q$). That is to say, any saturated set of 
elementary relations is either contained in $\Cal F$ or is all of $\Cal Q$. In terms 
of closed invariant sets, this means that any nonempty closed invariant set contains
$\Omega_{\Cal F}$. In particular $\Omega_{\Cal F}\subset \partial\Omega$, since $\partial\Omega$
is nonempty.

To complete the proof it is enough to show that $\Omega_{\Cal F}$ is nontrivial.
In fact, we shall show directly that $\partial\Omega\subset\Omega_{\Cal F}$. 
Suppose, by way of contradiction, that $\omega\in\Omega$ is a maximal directed set that
does not belong to $\Omega_{\Cal F}$. Then there is some relation $F\in\Cal F$ 
which $\omega$ fails. That is $\omega\cap tF = \emptyset$ for some $t\in \omega$.
By invariance of the set of maximal directed sets, we may replace 
$\omega$ by $t\inv \omega$, and thus assume that $\omega\cap F=\emptyset$. 
But then, for every $z\in\omega$, there is some $x_z\in F$ which 
has a common upper bound with $z$ (and with every element of $zP\inv$). 
Since $F$ is finite we may find a net of sets $zP\inv$ that
converges to $\omega$ and such that $x_z$ is constant, say $x_0$. But then
$\omega\vee x_0$ is a directed hereditary set which properly contains $\omega$, 
contradicting the maximality of $\omega$. 
\end{proof}

\begin{remark}
It is a consequence of the above proof that $\Cal F$ is a saturated
set of elementary relations. Also notice that one should not  expect the set of maximal directed
hereditary sets to be closed in $\Omega$. Indeed, in some cases
(cf. \proref{IrredCase}) $\partial\Omega$ is the entire Nica 
spectrum $\Omega$, which clearly contains elements that are not maximal. 
\end{remark}
\bigskip
 
Suppose now that $P$ has a finite set of lower bounds in the sense of Nica \cite{nica},
namely, assume there exists a finite set $S\subset P\setminus\{e\}$ such that 
\[
P\setminus\{ e\} = \bigcup_{s\in S} sP\,.
\]
In this case, Nica has shown that the projection 
$f_S=\prod_S(1-1_s)$ generates an ideal of the Toeplitz C*-algebra that is isomorphic 
to the compact operators on a separable, infinite dimensional, 
Hilbert space \cite[Proposition 6.3]{nica}. 
From the faithfulness theorem \cite[Theorem 3.7]{quasilat},
see also \cite[Theorem 24]{CL}, this ideal is minimal. Moreover,
the spectrum of this ideal may be characterised as follows:

\begin{definition}
The \emph{essential part of the Nica spectrum}, or more simply 
the \emph{essential spectrum}, is defined to be the spectrum $\Omega_\Cal E$ 
of the collection of elementary relations
\[
\Cal E := 
\{\ E\in\Cal Q\ :\ P\setminus \bigcup\limits_{x\in E} xP \text{ is a finite set }\}\,.
\]
Note that when $S$ is a finite set of lower bounds for $P$ then $S\in\Cal E$.
On the other hand, if $P$ has no finite set of lower bounds then $\Cal E$ has only the
trivial relations ($e\in E$ for all $E\in\Cal E$), in which case $\Omega_{\Cal E}=\Omega$.
\end{definition}
 
The next lemma gives a characterisation of the elements of $\Omega_{\Cal E}$
in the case where $P$ has a finite set of lower bounds.

\begin{lemma}\label{largest-invariant-set}
Suppose that $S$ is a finite set of lower bounds for $P$. Then
\[
\Omega_{\Cal E}=
\{\ \omega\in\Omega\ :\ \omega \text{ has no maximal element }\} = \Omega_{\{S\}} \,.
\]
Moreover, $\Omega_{\{S\}}$ is the (unique) largest proper closed invariant 
subset of $\Omega$ and the ideal of $\Cal T (G,P)$ generated by the function 
$f_S = \prod_{s\in S} (1 - 1_s)$ is the unique minimal ideal.
\end{lemma}
 
\begin{proof} These arguments are essentially from \cite{nica}, though adapted to include 
our definition of the essential spectrum $\Omega_{\Cal E}$.
Suppose that $\omega$ is a directed hereditary set which has no maximal element. 
Then, for every $t\in\omega$, the set $tP\cap\omega$ is infinite; if not then 
$\bigvee(tP\cap\omega)$ would be a maximal element of $\omega$. 
Now, let $E\in\Cal E$ and pick any $t\in\omega$. 
Since $tP\cap \omega$ is infinite it must intersect at least
one of the sets $txP$ for $x\in E$, and, since $\omega$ is hereditary, it follows
that $tP\cap \omega$ contains at least one element from $tE$.  
Thus $\omega$ satisfies the elementary relation $E$ at the value $t\in\omega$.
Since the above argument holds for any $t\in\omega$ and any $E\in\Cal E$ we have
that $\omega\in\Omega_{\Cal E}$.

On the other hand, if $\omega$ has a maximal element, $x$ say, 
then $\omega\cap xS=\emptyset$. In other words, $\omega$ fails 
to satisfy $S$ at the value $x \in \omega$. Thus $\omega\notin\Omega_{\{S\}}$. 

The above two arguments show that
$
\Omega_{\{S\}} \subseteq\{\ \omega\in\Omega\ :\ \omega \text{ has no maximal element }\}
\subseteq \Omega_{\Cal E}\,,
$
and equality follows from the observation that $S\in\Cal E$ and hence that 
$\Omega_{\Cal E}\subseteq \Omega_{\{S\}}$.
 
To see that  $\Omega_{\{S\}}$ is maximal, suppose
$\omega \notin \Omega_{\{S\}}$, then  
$\omega$ has a maximal element $x$, which satisfies $x\inv \omega \cap P = \{e\}$,
from which it follows easily that
 $\omega$ can satisfy only  those relations $H\in\Cal Q$ 
that contain the identity and are therefore everywhere satisfied. Thus if $X$ is not contained
in  $\Omega_{\{S\}}$ then $\Cal L(X)=\{ H\in\Cal Q : e\in H\}$ and the 
smallest closed invariant set $\Omega_{\Cal L(X)}$ containing $X$ is all of $\Omega$.  
\end{proof}

\begin{remark}
The set of elementary relations $\Cal E$ is,
in fact, saturated; to see this it suffices to verify that
if $H\notin \Cal E$ then $P\setminus \bigcup_{x\in H} xP$ contains an infinite ascending
chain. 
\end{remark}

\begin{remark}\label{partialomegaremark}
We have identified two sets of elementary relations, 
$\Cal E$ and $\Cal F$, generating extremal ideals in the lattice of 
invariant ideals of $\Cal T (G,P)$. In terms of closed invariant subsets 
of the Nica spectrum, we have the inclusions

\[   
\partial\Omega=\Omega_{\Cal F}\ \subseteq\  X\ \subseteq\  \Omega_{\Cal E} 
\]
for every proper closed invariant set $X\subset \Omega$. Moreover, on one side,
the boundary spectrum $\partial\Omega$ is always nonempty and, on the other side,
the essential  spectrum is a proper subset of $\Omega$ precisely 
when $P$ admits a finite set of lower bounds:
\[
\Omega_{\Cal E} =\left\{
\begin{aligned}
\Omega_{\{ S\}} \hsthree &\text{ if $P$ admits a finite set $S$ of lower bounds, }\\
\Omega\ \  \hsfour &\text{ otherwise. }
\end{aligned}\right.
\]
The reverse inclusions hold between the corresponding ideals of $\Cal T (G,P)$.
\end{remark}


\section{Controlled length functions and amenability}\label{SectAmenable}


 A sufficient condition for amenability of a quasi-lattice
order $(G,P)$ was established in \cite{quasilat}, it 
corresponds to amenability of the partial action of $G$ on
the Nica spectrum $\Omega$ from \cite{ELQ}. 
This was applied in \cite{CL} to prove amenability in the case where
$(G,P)$ is a right-angled Artin group. In this section we extend the argument of 
\cite{quasilat} to show that the same condition, namely, the existence of a
``controlled'' map of $G$ onto an amenable group, can be used to prove amenability 
of the canonical partial action of $(G,P)$ on \emph{all} closed invariant
subsets of the Nica spectrum $\Omega$. We begin by rephrasing the notion of 
generalised length function that lay at the heart of the argument in \cite{quasilat}.
 
\begin{definition}\label{controlledmapdef}
Suppose that $(G,P)$ and $(\Cal G,\Cal P)$ are quasi-lattice orders.
A {\em controlled map} is an order-preserving homomorphism 
$\varphi:(G,P)\to (\Cal G,\Cal P)$ such that
\begin{itemize}
\item[(C1)] the restriction $ \varphi_P:P\to\Cal P $ is finite-to-1, and
\item[(C2)] for all $x,y\in P$ satisfying $x\vee y\neq\infty$ 
we have $\varphi(x)\vee\varphi(y)=\varphi(x\vee y)$.
\end{itemize}
A controlled map $\varphi$ is most useful when it takes values in an amenable
(typically free abelian) group $\Cal G$. One should 
think of $\varphi$ as a type of generalised length function. 
\end{definition}

\noindent Note that condition (C1) implies the following property: 
\begin{itemize}
\item[(C3)] $\varphi^{-1}(e)\cap P=\{e\}$.
\end{itemize}
The reason is that if $\varphi(x)=e$ for some $x\in P\setminus\{e\}$ we would also have 
$\varphi(x^n)=e$ for all $n\in\N$, contradicting $\varphi_P$ finite-to-1. 
With this, one can easily deduce the following further properties of a controlled map
$\varphi$, for all $x,y\in P$:
\begin{itemize}
\item[(C4)] if $x\vee y\neq\infty$ and $\varphi(x)\leq\varphi(y)$ then $x\leq y$; and
\item[(C5)] if $x\vee y\neq\infty$ and $\varphi(x)=\varphi(y)$ then $x=y$.
\end{itemize}
Property (C5) follows easily from (C4), and property (C4) is deduced from (C2) and (C3),
by using the fact that $x\leq y$ if and only if $y^{-1}(x\vee y)=e$. 
We remark that only properties (C2) and (C5) were explicitly stated in  \cite{quasilat}.
However, property (C1) also holds for all examples considered in that paper and subsequently
in \cite{CL}; it is essential for the proof of \lemref{lemmAB} below. 

\begin{definition}
We say that a subset of $P$ is \emph{closed under taking least upper bounds}
(or simply \emph{closed under $\vee$}, or \emph{$\vee$-closed}) 
if it contains the least upper bound of any two
of its elements, whenever it exists. Define the \emph{$\vee$--closure} $F^{\vee}$ 
of $F \subset P$ to be the smallest $\vee$-closed subset of $P$ containing $F$.
Clearly, when $F$ is a finite subset of $P$, its $\vee$-closure may be written 
\[
F^{\vee} = \{\  \vee E\,:\, E\subseteq F \text{ has 
an upper bound  in } P\,\}\,,
\] 
which is again a finite set.
\end{definition}
   
\begin{lemma}\label{lemmAB}
Let $\varphi:(G,P)\to (\Cal G,\Cal P)$ be a controlled map and let $X$ be a closed 
invariant subset of the Nica spectrum $\Omega$ of $(G,P)$.
Suppose that $A$ is a $\vee$-closed finite subset of $\Cal P$
and that $s_0$ is a minimal element of $A$. 
Then either
\begin{itemize}
\item [{\bf (A)}] there exists $H\in \mathcal L(X)$ such that
$s_0\varphi(H)\subset A\setminus\{ s_0\}$ (in particular, $e\notin H$), or
\item [{\bf (B)}] there exists $\omega_0 \in X$ such that,
for all $x\in\varphi_P^{-1}(s_0)$ and all $y\in\varphi_P^{-1}(A\setminus\{s_0\})$,
we have $y\notin x\omega_0$. 
\end{itemize} 
\end{lemma}
\begin{proof}
Define the set
\[
H=\{\ x\inv(x\vee y)\ :\  x\in\varphi_P^{-1}(s_0)\,,
\ y\in\varphi_P\inv(A\setminus\{s_0\})\, \text{ and } x\vee y\neq\infty\ \}\,.
\]
Since the restriction $\varphi_P$ of the controlled map $\varphi$ to the 
semigroup $P$ is finite-to-1 (by property (C1) of Definition~\ref{controlledmapdef}) both $\varphi_P\inv(s_0)$ and
$\varphi_P\inv(A\setminus\{s_0\})$ are finite sets, and so $H$ is finite. 
Applying $\varphi$ to $H$ and left-multiplying by $s_0$ yields
\[
s_0\varphi(H)=\{\varphi(x\vee y) : x\in\varphi_P\inv(s_0)\,,
\  y\in\varphi_P\inv(A\setminus\{s_0\})\,,
\text{ and } x\vee y\neq\infty\}\,.
\]
By property (C2) 
elements of  $s_0\varphi(H)$ have the form
$\varphi(x\vee y)=\varphi(x)\vee\varphi(y) = s_0\vee \varphi(y)$ where
$\varphi(y)\in A\setminus\{ s_0\}$ and, since $A$ is $\vee$-closed and $s_0$ is
a minimal  element of $A$, we have that $s_0\varphi(H)\subset A\setminus\{ s_0\}$.

If $H\in\mathcal L(X)$, we have {\bf (A)}. Otherwise, there is 
some $\omega_0$ in $X$ that does not satisfy the elementary relation $H$. 
By invariance of $X$, we may choose $\omega_0$ such that $\omega_0 \cap H = \emptyset$.
Suppose that $x\in\varphi_P^{-1}(s_0)$ and $y\in\varphi_P^{-1}(A\setminus\{s_0\})$.
If $y \in x\omega_0$ then, because $x\omega_0$ is directed, we have $x\vee y\in x\omega_0$
(in particular $x\vee y\neq \infty$). But then $x^{-1}(x\vee y)\in \omega_0\cap H$, 
contradicting the choice of $\omega_0$. 
It follows that $ y \notin x\omega_0$, for all such $x,y$, and so {\bf (B)} holds.
\end{proof}
 
\begin{lemma}\label{PartialRep}
Let $(G,P)$ be a quasi-lattice ordered group and let $X$ denote a closed invariant 
subset of the Nica spectrum of $(G,P)$. Let $\{ \varepsilon_\omega : \omega\in X\}$
be the canonical orthonormal basis of the Hilbert space $\ell^2(X)$.
Then
\[
V_x(\varepsilon_\omega)= 
\begin{cases}\varepsilon_{x\omega} & \text{ if } x^{-1}\in\omega \\ 
0 & \text{ otherwise,} 
\end{cases}
\hsten\text{ for each }x\in G\,,
\]
defines a partial representation $V: G\to B(\ell^2(X))$
that satisfies the Nica relations $\Cal N$ together with the set of
elementary relations $\Cal L(X)$.
\end{lemma}

\begin{proof}
Let $u$ denote the partial representation of $G$ on $\ell^2(\X_G)$, given by
$u_x(\varepsilon_\omega) = \varepsilon_{x\omega}$ when $x\inv \in \omega$ and $0$ otherwise.
Because of the invariance of $X$, the subspace $\ell^2(X)$ of
$\ell^2(\X_G)$ is invariant under the partial representation $u$, so  
by restriction we obtain a partial representation
 $V= u\restr_{\ell^2(X)}$.
In order to see that $V$ satisfies the relations, notice first that,
for any $x\in G$, and any $\omega\in X$, we have 
\[
V_xV_x^*(\varepsilon_\omega)= 
\begin{cases}\varepsilon_\omega & \text{ if } x\in\omega \\ 
0 & \text{ otherwise.} 
\end{cases}
\]
It is now easily checked that, since $X\subset\Omega$, the Nica relations are satisfied 
(for the fact that $\omega$ is directed and hereditary implies that 
$x\vee y\in\omega$ if and only if both 
$x\in\omega$ and $y\in\omega$).

Now, if $H$ is a finite subset of $P$ representing an elementary relation, 
then $V$ satisfies $H$ if $\prod_{x\in H} (1-V_xV_x^*)=0$, 
equivalently if $\omega\cap H\neq\emptyset$ for every $\omega\in X$. 
But this is clearly the case whenever $H\in\Cal L(X)$.
\end{proof}

\begin{definition}
Let $(G,P)$ be a quasi-lattice ordered group and $X\subset\Omega$ a closed invariant subset
of the Nica spectrum.  By \cite[Theorem 4.4]{ELQ} and the fact that $X=\Omega_{\Cal L(X)}$,
the partial representation $V$, defined in Lemma \ref{PartialRep} above, determines a unique 
representation $\rho_V$ of $C(X)$ such that $(\rho_V,V)$ is a covariant pair.
We define the {\em spectral representation} of $C(X)\rtimes G$ to be 
\[
\lambda=\rho_V\times V: C(X)\rtimes G \to B(\ell^2(X))\,.
\]
In particular, $\lambda(i_X(x)i_X(y)^*) =V_xV_y^*$, where $i_X$ denotes the universal
partial representation of $G$ subject to the Nica relations and  
elementary relations $\Cal L(X)$. Recall that this partial representation $i_X$
generates the crossed product $C(X)\rtimes G$, \cite{ELQ}.
\end{definition}

\begin{proposition}\label{faithful}
Let $\varphi:(G,P)\to (\Cal G,\Cal P)$ be a controlled map, and suppose that
$X$ is a closed invariant subset of the Nica spectrum of $(G,P)$.
Then
\[
\Cal K := \overline{\rm span}\{ i_X(x)i_X(y)^* : x,y\in P \text{ with }
\varphi(x)=\varphi(y)\in\Cal P\}
\]
is a $C^*$-subalgebra of $C(X)\rtimes G$ on which $\lambda$ is faithful.
\end{proposition}

\begin{proof}
For each $\vee$-closed finite set $A\subset\Cal P$ define
\[
\Cal K_A : = \overline{\rm span} \{ i_X(x)i_X(y)^* : x,y\in P \text{ with } 
\varphi(x)=\varphi(y)\in A\}\,,
\]
and write $\Cal K_s$ for $K_{\{ s\}}$. Clearly 
$\{\Cal K_A : A\subset\Cal P\text{ finite, $\vee$-closed}\}$, 
with the $A$'s directed by inclusion,
is an inductive system with limit $\Cal K = \overline{\cup_A \Cal K_A}$.
We need to check that each $\Cal K_A$ is closed under multiplication. 
Note first that, by the Nica relations, 
\begin{equation}
i_X(x) i_X(y)^* i_X(z) i_X(w)^* = 
\begin{cases}
i_X(x y\inv (y\vee z)) i_X(w z\inv (y\vee z))^* &\text{ if } y\vee z \neq\infty \\ 
0 & \text{ otherwise.} 
\end{cases}
\end{equation}
If  $\varphi(x)=\varphi(y)\in A$, $\varphi(z)=\varphi(w)\in A$ 
and $y \vee z \neq\infty$, then  property (C2) 
of the controlled map $\varphi$,
 and the fact that $A$ is $\vee$-closed,
implies that $\varphi(x y\inv (y\vee z))= 
\varphi(w z\inv (y\vee z))=\varphi(y\vee z)=\varphi(y)\vee\varphi(z)\in A$.
It follows that  each $\Cal K_A$ is a $C^*$-subalgebra of $C(X)\rtimes G$,
and so is their limit $\Cal K$,
cf. \cite[Lemma 4.1]{quasilat}. 
To conclude that $ \lambda $ is faithful on $\Cal K$, it suffices to show
that $\lambda$ is faithful on $\Cal K_A$ for each $\vee$-closed finite $A$,
 and then apply \cite[Lemma 1.3]{alnr}.

Fix some $\omega_0\in X$. For each $s\in\Cal P$, we define
\[
H_s^{(\omega_0)} =  \overline{\rm span}
\{\,\varepsilon_{z\omega_0}\in\ell^2(X)\ :\ z\in P\,,\  \varphi(z)=s\,\}\,.
\]
Let $\pi_s^{(\omega_0)}$ denote the projection onto $H_s^{(\omega_0)}$.
Although not strictly essential to the proof, cf. \cite{quasilat}, it is interesting to note that
t $H_s^{(\omega_0)}$ is a finite dimensional subspace of $\ell^2(X)$, 
by property (C1).

Suppose we are given $x,y\in P$ such that $\varphi(x)=\varphi(y)$. 
 Using  property (C5)
of the controlled map  $\varphi$ and the fact that $\omega_0$ is directed,
we have that  $x\omega_0=y\omega_0$ if and only if $x=y$.
Consequently, for $s\in \Cal P$, the orthonormal basis 
$\{\varep_{z\omega_0} : z\in\varphi_P^{-1}(s)\}$
for $H_s^{(\omega_0)}$ is in bijection with $\varphi_P^{-1}(s)$.
Moreover, if $x,y,z\in \varphi_P^{-1}(s)$, we have 
\begin{equation}\label{action}
\lambda(i_X(x)i_X(y)^*)(\varepsilon_{z\omega_0}) =V_xV_y^*(\varepsilon_{z\omega_0})=
\begin{cases}\varepsilon_{x\omega_0} & \text{ if } y = z\\ 
0 & \text{ otherwise.} 
\end{cases}
\end{equation}
Thus, given $s\in\Cal P$ and any $x,y\in P$ such that $\varphi(x)=\varphi(y)=s$,
the operator $\lambda(i_X(x)i_X(y)^*)$ 
restricts to the rank--one operator
$\<\ \cdot\ ,\varepsilon_{y\omega_0}\>\varepsilon_{x\omega_0}$ on the finite 
dimensional Hilbert space $H_s^{(\omega_0)}$.  
It follows that $\lambda(\cdot) \restr H_s^{(\omega_0)}$ is an isomorphism between 
$\Cal K_s$ and the algebra of operators on 
$H_s^{(\omega_0)}$.
So $\lambda$ is isometric on $\Cal K_s$ for each $s\in\Cal P$, and in particular, 
if $T_s\in\Cal K_s$ then $\|\lambda(T_s)\pi_s^{(\omega_0)}\|=\|T_s\|$. 
This last statement will be needed below.
\smallskip

We now show by induction on the size of $A$ that $\lambda$ is faithful on $\Cal K_A$
for each $\vee$-closed finite subset $A\subset\Cal P$. Let $A$ be such a subset of $\Cal P$,
and let $s_0$ be a minimal element of $A$.

There are two cases to be considered, according to \lemref{lemmAB}:
in case {\bf (A)}  we will show that $\Cal K_A=\Cal K_{A\setminus\{ s_0\}}$
and in case {\bf (B)} we will show that if $T\in\Cal K_A$ and $\lambda(T)=0$
then $T\in\Cal K_{A\setminus\{ s_0\}}$. 
Since $A\setminus\{s_0\}$ is also finite and closed under $\vee$ because $s_0$
is minimal, the result will then follow by induction.
\medskip

\noindent{\bf Case (A):} 
There exists $H\in\Cal L(X)$ such that $s_0\varphi(H)\subset A\setminus\{ s_0\}$. 
Note that the partial representation $x\mapsto i_X(x)$ satisfies all
the relations $\Cal L(X)$ and $\mathcal N$. In particular, we have 
\begin{equation}\label{relH}
\prod_{u\in H} (1-i_X(u)i_X(u)^*) =0 \,.
\end{equation}
Using the Nica relations, Equation \eqref{relH} may be re-expressed as
\begin{equation}\label{relH2}
1=\sum_{z\in H^\vee} n_zi_X(z)i_X(z)^*\ , \qquad \text{ where }n_z\in\Z\,,
\end{equation}
and where the sum is taken over the $\vee$-closure $H^\vee$,
which is a finite set.  
By property (C2) of the controlled map $\varphi$, 
we have $s_0\varphi(H^\vee)= (s_0\varphi(H))^\vee$ and, since 
$A\setminus\{ s_0\}$ is $\vee$-closed, it follows that
$s_0\varphi(H^\vee)\subset A\setminus\{s_0\}$. 
Now, for any $x,y\in\varphi_P^{-1}(s_0)$, Equation \eqref{relH2} gives
 $i_X(x)i_X(y)^* =\sum\limits_{z\in H^\vee} n_z i_X(xz)i_X(yz)^*$ which,
by the previous statement, is an
element of $\Cal K_{A\setminus\{ s_0\}}$. Therefore 
$\Cal K_A = \Cal K_{A\setminus\{ s_0\}}$.
\medskip

\noindent{\bf Case (B):}
There exists $\omega_0\in X$ such that $y\notin z \omega_0$ 
whenever $z\in\varphi_P^{-1}(s_0)$ and $y\in\varphi_P^{-1}(A\setminus\{s_0\})$.
Therefore, given $x,y,z\in P$ such that $\varphi(x)=\varphi(y)\in A\setminus\{s_0\}$
and $\varphi(z)=s_0$, we have that $V_xV_y^*(\varep_{z\omega_0})=0$
(since $y\notin z\omega_0$). It follows 
that $\lambda(\Cal K_{A\setminus\{ s_0\}})\pi_{s_0}^{(\omega_0)}=0$.
Suppose now that $T\in\Cal K_A$ satisfies $\lambda(T)=0$, and write 
\[
T=\lim_n\sum_{s\in A} T_{n,s}\quad\text{ where }
T_{n,s}\in\Cal K_s \text{ for } s\in A\,.
\]
Then 
\[
\lambda(\sum_{s\in A} T_{n,s})\pi_{s_0}^{(\omega_0)}
= \lambda(T_{n,s_0})\pi_{s_0}^{(\omega_0)} \longrightarrow 0\,,\hsthree
\text{ as } n \longrightarrow\infty\,.
\]
But since $\|\lambda(T_{n,s_0})\pi_{s_0}^{(\omega_0)}\|=\|T_{n,s_0}\|$
we have $\|T_{n,s_0}\|\to 0$. Therefore $T\in\Cal K_{A\setminus\{ s_0\}}$.
\end{proof}

\begin{theorem}\label{amenable}
Let $\varphi:(G,P)\to (\Cal G,\Cal P)$ be a controlled map and let $X$ be a closed 
invariant subset of the Nica spectrum $\Omega$ of $(G,P)$.
If $\Cal G$ is an amenable group then the conditional 
expectation $\Phi: C(X)\rtimes G\to C(X)$ is faithful on positive elements
and so the canonical partial action of $G$ restricted to $C(X)$ is amenable.
\end{theorem}

\begin{proof} 
The proof uses an adaptation of the argument 
from \cite[Proposition 4.2]{quasilat} to the present situation.
Recall that the elements of the form $i_X(x)i_X(y)^*$ span a dense 
*-subalgebra of the crossed product $C(X)\rtimes G$ and that 
on these elements, the conditional expectation $\Phi$ is given 
by $\Phi(i_X(x)i_X(y)^*)= i_X(x)i_X(x)^*$ if $x = y$ and $\Phi(i_X(x)i_X(y)^*)= 0$ if $x\neq y$.  
The controlled map $\varphi:G \to \Cal G$ induces a coaction of $\Cal G$ on 
 $C(X)\rtimes G$, determined  by 
$i_X(x) i_X(y)^* \mapsto i_X(x) i_X(y)^* \otimes \delta_{\varphi(xy\inv)}$, 
whose associated conditional expectation is 
determined by 
\[
\Phi_{\Cal G} (i_X(x) i_X(y)^*)= 
\begin{cases} i_X(x) i_X(y)^* & \text{ if } \varphi(x) = \varphi(y) \\ 
0 & \text{ otherwise.} 
\end{cases}
\] 
The fixed point algebra of this coaction is the C*-algebra 
$\Cal K$ described in \proref{faithful},
and $\Phi_{\Cal G}: C(X)\rtimes G\to \Cal K$ 
is faithful on positive elements because $\Cal G$ is amenable. 
Since $\Phi = \Phi \circ \Phi_{\Cal G}$, 
in order to conclude that the partial action 
of $G$ on $X$ is amenable, we need
to prove that $\Phi: \Cal K \to C(X)$ is faithful on positive elements.
This is done using the faithful representation $\lambda$ 
of $\Cal K$ on $\ell^2(X)$ from \proref{faithful}.

For each $\omega\in X$ let $P_\omega$ denote the rank-one projection onto 
the basis vector $\varepsilon_\omega$ of $\ell^2(X)$, and 
 consider the diagonal map $\Delta$ on  $ B(\ell^2(X))$ defined by 
$\Delta (U) = \lim_I(\sum_{\omega\in I} P_\omega U P_\omega)$, 
where the weak limit is taken over the finite subsets of 
the canonical orthonormal basis of $\ell^2(X)$ directed by inclusion. 
 Then, for a positive element $U=b^*b$, 
we have $\Delta(U) = 0$ only if $U = 0$. In particular, if 
$T\in \Cal K_+$, then $\Delta(\lambda(T)) = 0$ only if $\lam(T)=0$ and so,
by faithfulness of $\lam$, only if $T=0$. It now suffices to show that the 
compositions $\Delta\circ\lambda$ and $\lambda\circ\Phi$ agree on $\Cal K$, for then
the fact that $\Delta\circ\lam$ is faithful on $\Cal K_+$ implies that $\Phi$ is also 
faithful on positive elements of the range of $\Phi_{\Cal G}$.

From the definition of $\Phi$, it is clear that 
$\lambda\circ\Phi(i_X(x) i_X(y)^*)$ is equal to $ V_x V_x^*$ when $x=y$ and is
$0$ otherwise. Next recall that for any $\omega\in X$ and $s\in\Cal P$ 
the set $\varphi_P\inv(s)\cap\omega$ contains at most one element
(because $\omega$ is directed and $\varphi$ satisfies property (C5) of a controlled map).
Thus, given $x,y\in P$ such that $\varphi(x)=\varphi(y)=s$, 
we have, by \eqref{action}, that
\[
P_\omega V_x V_y^* P_\omega = 
\begin{cases}
P_\omega &\text{ if }x = y\in\omega\\
0 &\text{ otherwise.}
\end{cases}
\]
 Therefore $\Delta\circ\lambda (i_X(x) i_X(y)^*) = \Delta (V_x V_y^*)$ is given by 
\[
\lim_I (\sum_{\omega\in I} P_\omega V_x V_y^* P_\omega)\ =\ 
\begin{cases}
\ \lim_I (\sum_{\omega\in I} 1_x(\omega) P_\omega) = V_xV_x^* &\text{ if }x = y\\
\ 0 &\text{ otherwise.}
\end{cases}
\]
This proves that $\lambda\circ\Phi = \Delta\circ\lambda$ on $\Cal K$, as required. 
\end{proof}

\begin{corollary}\label{coramenable}
Assume now that $\Gamma$ is a simplicial graph with vertex set $S$. Let
$A_\Gamma$ be the corresponding right-angled Artin group and let $\Omega$ the 
Nica spectrum for the quasi-lattice order $(A_\Ga,A_\Ga^+)$. 
Then the abelianisation map  $\phi: (A_\Gamma, A_\Gamma^+) \to (\Z,\N)^{|S|}$ 
is a controlled map, and the partial action of $A_\Gamma$ is amenable 
when restricted to any closed invariant subset of $\Omega$.
\end{corollary}

\begin{proof}
Property (C2) of a controlled map 
was proved in \cite[Proposition 19]{CL}. To check that property (C1) holds
it suffices to consider the length homomorphism 
$\ell : (A_\Ga,A_\Ga^+)\to (\Z,\N)$ such that $\ell(s)=1$ for each
positive generator $s\in S$. 
If $\phi(x) = s$ and $\ell(x) = n$, then each 
positive element $y \in A_\Ga^+$ such that $\phi(y) = s$
also has length $n$ and moreover, is a product of
the positive generators that appear in $x$, with the same repetitions,
so there are at most $n!$ such $y$'s, proving that $\phi$ satisfies (C1).
Amenability of the restricted partial action now follows from \thmref{amenable}.
\end{proof}

%
\section{Topological freeness on $\partial\Omega$ and simplicity of the boundary quotient}
\label{SectBoundary}
%

We are interested in deciding when the boundary quotient $C(\partial\Omega)\rtimes G$ of
the Toeplitz algebra of a quasi-lattice ordered group is purely infinite and simple.
To this purpose, we first give a theorem
for general quasi-lattice ordered groups, unifying and strengthening similar
results from \cite{quasilat,ELQ,CL}. One of the 
key assumptions of the theorem is the amenability discussed in the preceding section;
the other one is topological freeness, and in the remainder of the section 
 we develop the necessary tools to decide which right-angled Artin groups have this 
property. 

\begin{theorem}\label{PureSimple}
Let $(G,P)$ denote a quasi-lattice ordered group. If the partial action
of $G$ on $\partial\Omega$ is both amenable and topologically free then 
the boundary quotient $C(\partial\Omega)\rtimes G$ is purely infinite and simple.
\end{theorem}

\begin{proof}
It follows  from Proposition 4.3 of \cite{purelinf}
that the canonical partial action of $G$ on $C(\partial\Omega)$ is minimal.
By  Corollary 2.9 of \cite{ELQ} and topological freeness,
the reduced crossed product is simple. By amenability of the action,
the reduced crossed product is isomorphic to the full crossed product
$C(\partial\Omega)\rtimes G$. In order to show that the algebra is 
purely infinite we adapt a familiar argument from \cite{cun}, along the lines of 
the proof of \cite[Theorem 5]{lac-spi}.

  Let $x\neq 0$ be an element of $C(\partial \Omega)\rtimes G$. 
Since the canonical conditional expectation $\Phi$ is faithful, $\Phi(x^* x) \neq 0$ 
and we may define $a = x^* x / \|\Phi(x^* x)\|$, so that 
$\|\Phi(a)\| = 1 $. We may find a finite sum 
$b = \sum_{t\in F}  b_t\cdot \delta_t  \in  C_c(G, C(\partial \Omega))$, where each
 $b_t$ is in the domain $D_t=C_0(U_t)$ of the partial 
isomorphism $\alpha_t$,  such that $b \geq 0$ and  $\|b-a\|<  \frac{1}{4}$. 
Since $\Phi$ is positive and contractive, $\Phi(b) = b_e$ is a positive element of 
$C(\partial \Omega)$ with norm greater than  $3/4$.

Now let $U$ be a nonempty open subset of $\partial \Omega$ 
on which the function $b_e$ takes values greater than $3/4$. 
By topological freeness there exists a nonempty open subset $U_0 \subset U$
such that for every $t\in F$ we have that $tU_0 \cap U_0 = \emptyset$.
Notice that this may happen either because $U_0$ does not meet the domain 
$U_{t\inv}$ of $\theta_t$ or else because  
$tU_0 = \theta_t(U_0\cap U_{t\inv})$ is nonempty but disjoint from $U_0$.

Let $V(t,tH)$ be a basic open set of $\Omega$ such that $V(t,tH)\cap\partial\Omega$ 
is nonempty and contained in $U_0$. Without loss
of generality we may also suppose that $t\in P$. Indeed, $V(t,tH) \neq \emptyset$ 
implies that $t\vee e\neq\infty$ (since $t\in\omega$ for $\omega\in V(t,tH)$),
and replacing $t$ by  $u:= t\vee e \in P$ and 
$H$ by  $K:=\{ u^{-1}(th\vee e):h\in H \text{ and }th\vee e\neq\infty\}$ we have that $V(t,tH)=V(u,uK)$.

Since $V(t,tH)\cap\partial\Omega\neq\emptyset$ we must have $H\notin\Cal F$, 
and so we may find $z\in P$ such that $z\vee h=\infty$ for all $h\in H$. 
Setting $y=tz$ (an element of $P$ since $t\in P$), we conclude that if
$\omega\in\partial\Omega$ contains $y$, then  $ \omega\in V(t,tH)\subset U_0$.  
In other words, the characteristic function of $U_0$ dominates the projection 
$1_y$ with $y\in P$. 
We note in passing that this property plays the role of the boundary action property 
needed in the proof of \cite[Theorem 5]{lac-spi}.
It now follows that $b_e \geq \frac{3}{4} 1_y $ and thus 
\[
V_y^* b V_y = V_y^* b_e V_y 
\geq V_y^* ({\textstyle\frac{3}{4}}V_y V_y^*)V_y = {\textstyle\frac{3}{4}} I
\]
Since $\| V_y^* a V_y -  V_y^* b V_y\| < \frac{1}{4}$, 
the element $V_y^* a V_y$ is invertible, and setting 
\[
A := \|\Phi(x^* x)\|\inv (V_y^* a V_y)\inv V_y^* x^*
\quad \text{ and } \quad B:= V_y,
\] 
and since $a = x^* x / \|\Phi(x^* x)\|$ we obtain
\[
AxB  = \|\Phi(x^* x)\|\inv (V_y^* a V_y)\inv V_y^* x^* x V_y = I.
\]
\end{proof}

In order to characterise topological freeness for quasi-lattice groups, 
we need to introduce the following notion.
\begin{definition}
We shall say that a submonoid $P'$ of a monoid $P$ is a \emph{full submonoid}
if whenever $x\in P'$ and $x=ab$ for $a,b\in P$ we have both $a$ and $b\in P'$.
\end{definition}
 
 \begin{lemma}\label{coreLemma}
Let $(G,P)$ be a quasi-lattice ordered group and define  
\[
P_0 :=\{\ x\in P \ :\ x\vee y\neq\infty \text{ for all } y\in P\ \}\,.  
\] 
Then $P_0$ is a directed full
submonoid of $P$ and generates a subgroup of $G$ of the form $G_0=P_0P_0\inv$. 
Consequently, $(G_0,P_0)$ is a lattice ordered group and the inclusion
into $(G,P)$ is an order preserving map. \end{lemma}

\begin{proof}
If $x,y\in P_0$ then $x\vee y\neq\infty$ and, for all $z\in P$, we have 
$y\vee z\neq\infty$, and $x\vee (y\vee z)\neq\infty$. Therefore
$x\vee y\in P_0$. So $P_0$ is directed. Similarly, 
$xy\vee z=x(y\vee(x\inv(x\vee z)))\neq\infty$, for all $z\in P$, and so $P_0$ is a submonoid.
A similar argument also shows that $P_0$ is a full submonoid of $P$, namely,
if $x=ab\in P_0$, for $a,b\in P$, then both $a$ and $b\in P_0$. 

To see that $G_0: =P_0P_0\inv$ is indeed a subgroup, observe that if $x,y\in P_0$ then 
$x\vee y\in P_0$ and by fullness
we may write $x\vee y=xa=yb$ for $a,b\in P_0$. But then $x\inv y=ab\inv$. Using this, 
and the fact that $P_0$ is a submonoid, every product of elements of $P_0$ 
and their inverses may be rearranged into the form $uv\inv$ for $u,v\in P_0$.
The remaining statements follow easily:
 in particular, if $z\in G_0$ then $z\leq a$
for some $a\in P_0$ which implies that $x\vee e\neq\infty$ and, by fullness of $P_0$,
that $x\vee e\in P_0$. It follows that, for all $x,y\in G_0$, 
we have $x\vee y=x(x\inv y\vee e)\in G_0$. 
\end{proof}

\begin{definition}
Let $(G,P)$ be a quasi-lattice ordered group. We shall refer to the sublattice 
$(G_0,P_0)$ introduced in Lemma \ref{coreLemma} as the \emph{core} of $(G,P)$.
\end{definition}

\begin{proposition}\label{TopFreeCore}
Let $(G,P)$ be a quasi-lattice ordered group with Nica spectrum $\Omega$ and 
core $(G_0,P_0)$.  Then the canonical partial action
of $G$ on its boundary spectrum $\partial\Omega$ is topologically free if and only
if the action of the core $G_0$ on $\partial\Omega$ is topologically free.
\end{proposition}

\begin{proof}
Clearly the action of $G_0$ is topologically free if the action of $G$ is.
In order to prove the converse, suppose that $\Fix_{\partial\Omega}(t)$ has nonempty interior
for some $t\in G\setminus \{ e\}$. We shall show that there exists some 
$t_0\in G_0\setminus\{e\}$ for which $\Fix_{\partial\Omega}(t_0)$ has
nonempty interior. 
Let $U\subset\partial\Omega$ be an open set pointwise fixed by $t$.
Since $\partial\Omega$ is the closure of the set of
maximal hereditary sets, we may find some $\omega\in U$ that is maximal.
We now claim that there is some $b\in\omega$ such that $t_0=b\inv tb$ lies in $G_0$.
Since $b\in\omega$, the domain $U_b$ of the partial 
homeomorphism $b\inv$ contains $\omega$, so intersects $U$ nontrivially. 
It follows that $b\inv(U_b\cap U)$ is a nontrivial open set fixed
pointwise by $t_0$, as required.     

To prove the claim we suppose, by way of contradiction, that $b\inv tb\notin G_0$ 
for all $b\in \omega$. Fix $b\in\omega$. Since $\omega$ is fixed by $t$ we also have 
$tb\in\omega$ and $t\inv b\in\omega$. Let $a:=b\inv(b\vee tb)$ and $c:=b\inv(t\inv b\vee b)$.
Then $ac\inv=b\inv t b$. Moreover, since $\omega$ is directed, both $ba$ and $bc$ lie
in $\omega$. Since $ac\inv\notin G_0$ we may suppose that one of $a$ or $c$ is not 
an element of $P_0$. If $a\notin P_0$ then there exists $z\in P$
such that $a\vee z=\infty$. Setting $\omega_b:=bz\omega\in\partial\Omega$ we have that 
$b\in \omega_b$ but $\omega_b\notin\Fix(t)$, for otherwise we would have $b,tb\in\omega_b$
and, by directedness, $ba=b\vee tb\in\omega_b$ and
$ba\vee bz \in\omega_b$ contradicting $a\vee z=\infty$.
Similarly, if $c\notin P_0$ then we can find $\omega_b\in\partial \Omega$ such
$b\in \omega_b$ but $\omega_b\notin\Fix(t)$. Since we can do this for every $b\in\omega$
we have  a net $\{\omega_b\}_{b\in\omega}$ of elements
of $\partial\Omega$, none of which are fixed by $t$.  Since
 $bP^{-1}\subset\omega_b$ for each $b\in\omega$, any limit point
of this net must contain $\omega$, and by maximality of $\omega$
we deduce that the net converges (uniquely) to $\omega$. 
This, however, contradicts the fact that $\omega$ lies in the interior
of $\Fix_{\partial\Omega}(t)$.
\end{proof}

Again, our main application is to right-angled Artin groups.
\begin{remark}
Note that the centre $A_0$ of a right-angled Artin group $A$
acts trivially on $\partial\Omega$.
Thus the canonical partial action of $A$ on its boundary spectrum $\partial\Omega$ fails 
to be topologically free whenever $A_0$ is nontrivial. 
For example, the boundary action of $(\Z,\N)$ is clearly not 
topologically free, as in this case $\partial\Omega$ is just a single point. 
It is easily seen that the core $(A_0,A_0^+)$ of a right-angled 
Artin group $(A,A^+)$ is just its maximal abelian direct factor $(\Z,\N)^n$;
in fact $A_0$ is the center of $A$.
 See Example~\ref{gammaopexample} below for a simple characterisation in terms of 
 the simplicial graph defining $A$.
 \end{remark}

\begin{corollary}\label{TopFreeArtin}
Suppose that $(A,A^+)$ is a right-angled Artin group that 
has no direct factor $(\Z,\N)$. Then $C(\partial\Omega)\rtimes A$
is purely infinite and simple.
\end{corollary}

\begin{proof}

If $(A,A^+)$ has no direct factor $(\Z,\N)$ then its centre (or core) is trivial
and, by \proref{TopFreeCore}, the canonical partial action of 
$A$ is topologically free on $\partial\Omega$. Since the abelianisation
of $(A,A^+)$ is a controlled map to an amenable group the canonical partial 
action is also amenable on $\partial \Omega$ by Theorem \ref{amenable}.
The result now follows by Theorem \ref{PureSimple}.
\end{proof}


\section{Irreducibility of $(G,P)$ and maximality of $\partial\Omega$}\label{SectIrredCase}


Let $P$ be a monoid. An \emph{atom} of $P$ is any nontrivial element $z\in P$ such that
whenever $z=xy$ for $x,y\in P$ then either $x=e$ or $y=e$. 
Note that when $P$ is the positive cone of a partially ordered group $(G,P)$ 
the atoms of $P$ are precisely the minimal elements of $P\setminus\{ e\}$ 
with respect to the partial order on $G$. By a \emph{set of lower bounds} for $P$
we mean a set $S\subset P$ such that $P\setminus\{ e\}=\bigcup_{s\in S} sP$ 
(cf. the definition of a \emph{finite} set of lower bounds 
in the sense of Nica, given in \secref{SectNica}).
Clearly, any set of lower bounds for $P$ will contain the set of all atoms of $P$.

In the following we denote by $(\Z,\N)^\infty$ the direct sum of countably infinitely many 
copies of the ordered integers $(\Z,\N)$. Note that we shall write the group operation  in 
$(\Z,\N)^\infty$ additively. We shall be interested here in quasi-lattice ordered groups 
which admit a controlled map to some free abelian group of finite or countably infinite rank,
as is the case for any right-angled Artin group $(A_\Gamma,A_\Gamma^+)$, by \corref{coramenable}.

\begin{lemma}\label{UsefulLemma}
 Let $(G,P)$ be a quasi-lattice ordered group which admits a controlled 
map 
\[
\phi:(G,P)\to (\Z,\N)^\infty\,.
\]
Then
\begin{itemize}
\item[(i)] The set of atoms of $P$ is a set of lower bounds for $P$.
\item[(ii)] There exists a homomorphism $\ell:P\to\N$ with $\ell^{-1}(0)=\{e\}$ 
such that whenever  $x\vee y\neq\infty$ (for $x,y\in P$) one has the inequality
$\ell(x\vee y)\leq \ell(x)+\ell(y)$.
\item[(iii)] If $a\in P$ is an atom and $e\leq x\leq a^k$, for some $k\in\N$,
then  $x=a^m$ for some $m\leq k$.
\end{itemize}

\end{lemma}

\begin{proof}
We first prove part (ii).
Let $\mu:(\Z,\N)^\infty\to(\Z,\N)$ denote the unique order homomorphism that 
maps each factor isomorphically onto $(\Z,\N)$ ,
thus $\mu((m_\lambda)_\lambda)=\sum_\lambda m_\lambda$ and $\mu$ satisfies part (ii).
We define the homomorphism $\ell:(G,P)\to (\Z,\N)$ by the composition
$\ell=\mu\circ\phi$. That $\ell$ satisfies part (ii) 
now follows because the map $\phi$ satisfies
properties (C2) and (C3) of a controlled map, Definition~\ref{controlledmapdef}.

In order to prove part (iii),  suppose $a\in P$ is an atom and $x\in P$ such that
$e \leq x \leq a^k$ for some $k\in \N$.
Notice that (iii) clearly holds in the case $k=1$. Assume now $k>1$.
Clearly  $a\vee x\leq a^k$, and from writing $a\vee x=ac$,
it follows that $c\leq a^{k-1}$. By induction, we may suppose that $c=a^r$ for some $r\leq k-1$, and
hence $a\vee x=a^{r+1}$. Let now $A=\phi(a)$ and $B=\phi(x)$, so
by property (C2) of Definition \ref{controlledmapdef},
 $A\vee B=\phi(a\vee x)=\phi(a^{r+1})=(r+1)A$.
 Recall that  when $A=(a_\lambda)_\lambda$ and
$B=(b_\lambda)_\lambda$,  then 
$A\vee B=(\max\{a_\lambda, b_\lambda\})_\lambda$. Since $A\vee B=(r+1) A$, we have that 
either $A \leq B = (r+1) A$ or else $B\leq A$, $B\neq A$ and $r = 0$.
By property (C4) of the controlled map $\phi$, and since $a$ is an atom, 
this implies that $x = a^{r+1}$ with $r+1§ \leq k$  or $x = e$, finishing the  proof of (iii).

For the proof of part (i), it suffices to show that every nontrivial element of
$P$ is divisible by an atom. Let $x\in P\setminus\{ e\}$. If $x$ is not itself
an atom we may write $x=yz$ where $y,z\in P\setminus\{e\}$. 
However, since the homomorphism $\ell$ is additive and nondegenerate 
($\ell^{-1}(0)=\{e\}$), we have $\ell(y)<\ell(x)$ and, by induction,
$a\leq y\leq x$ for some atom $a$. 
\end{proof}

\begin{remark}
Statement (i) of Lemma \ref{UsefulLemma} may be strengthened considerably: the monoid $P$ 
is actually generated by its atom set, and is \emph{atomic} in the sense of \cite{DP}, namely
for each $x\in P$ there is an upper bound on the length of any expression for $x$
as a product of atoms. However we will not need to make use of this extra information here.
\end{remark}


\begin{definition}
Suppose that $(G,P)$ is a quasi-lattice ordered group 
whose atom set $S$ is a set of lower bounds for $P$. 
We define the \emph{$\vee$-graph} of $(G,P)$ to be the graph with
vertex set $S$ and edges $\{ a,b\}$ whenever $a\vee b=\infty$.
We say that $(G,P)$ is \emph{graph-irreducible} if its $\vee$-graph is connected.  
\end{definition}

Note that if $(G,P)=(G_1,P_1)\oplus(G_2,P_2)$ is a proper direct sum of 
quasi-lattice ordered groups,  and if $S_i$ denotes the set of atoms of $P_i$ for $i=1,2$,
then the set of atoms of $P$ is just the union $S=S_1\cup S_2$. Moreover, $S$ is a set
of lower bounds for $P$ if and only if $S_i$ is a set of lower bounds
for $P_i$, for each $i=1,2$. 
In this case the $\vee$-graph 
for $(G,P)$ is just the disjoint union of the $\vee$-graphs for 
$(G_1,P_1)$ and $(G_2,P_2)$. Thus, graph-irreducibility is a 
stronger property than the more usual notion of irreducibility with respect to direct
sums.

Since all proper closed invariant subsets of the Nica spectrum of a quasi-lattice order
lie between the boundary spectrum $\partial\Omega$ and the essential 
spectrum $\Omega_{\Cal E}$, they are most easily classified when 
$\partial\Omega = \Omega_{\Cal E}$.
The following proposition describes circumstances under which this actually occurs.

\begin{proposition}\label{IrredCase}
Let $(G,P)$ be a quasi-lattice ordered group that admits a controlled 
map 
\[
\phi:(G,P)\to (\Z,\N)^\infty\,.
\]
If $(G,P)$ is graph-irreducible, then the boundary spectrum $ \partial\Omega$ 
coincides with the essential spectrum $\Omega_{\Cal E}\,$. 
If, in addition, $P$ has infinitely many atoms, 
then $\partial\Omega = \Omega$,  in which case the boundary quotient 
coincides with the universal Toeplitz C*-algebra $C^*(G,P)$. 
\end{proposition}

\begin{proof}
{}From Remark~\ref{partialomegaremark},  $\Omega_{\Cal E} =\Omega$
unless $P$ has a finite set of lower bounds, and clearly any set 
of lower bounds for $P$ contains all the atoms of $P$.
Therefore the second assertion of the proposition 
follows easily from the first.

By \lemref{minimalityofboundary}, we have $\partial\Omega \subset  \Omega_{\Cal E}\,$,  
so to prove the first assertion  it suffices to show that  
$\partial\Omega \supset  \Omega_{\Cal E}\,$, that is, $\Cal F\subset\Cal L(\Omega_\Cal E)$.
We therefore take a finite subset $F $ of $P$ such that  $F\notin \Cal L(\Omega_\Cal E)$ 
and we shall show that $F\notin\Cal F$.
We first introduce some notation. Given $x,y\in P$ we set 
\[
y\bsl x :=\left\{
\begin{aligned}
y^{-1}(y\vee x) \quad &\text{ if $y\vee x\neq\infty$, }\\
e \hsseven &\text{ otherwise. }
\end{aligned}\right.
\]
Let $P_F:=P\setminus\cup\{xP:x\in F\}$ and notice that if $x\in F$ 
and $y\in P_F$ then $y\bsl x = e$ if and only if $x\vee y = \infty$.
Let $S$ denote the set of all atoms of $P$ and note that, by Lemma \ref{UsefulLemma}(i),
this is a set of lower bounds for $P$.
For each $y\in P_F$, we also define the following subset of $S$:
\[
C(y):=\{ s\in S : s\leq\, y\bsl x \text{ for some } x\in F \} \,.
\]
Note that, since every nontrivial element of $P$ is bounded below by at least one atom,
the set $C(y)$ is non-empty unless $y\bsl x=e$ for all $x\in F$. 
Our objective is to find an element $y\in P_F$ for which $C(y)=\emptyset$; 
this will imply that $y\vee x=\infty$ for all $x\in F$ and hence that $F\notin\Cal F$. 
We break the proof into three steps.
\medskip

\paragraph{\bf Step 1.} \emph{The assumption $F\notin\Cal L(\Omega_\Cal E)$
implies that the set $P_F$ contains a strictly increasing infinite 
sequence of elements:
\[
y_1 < y_2 < \dots < y_i < y_{i+1} \dots\ ,\hsthree i\in\N\,.
\]}
 
In the case that $S$ is infinite, this is a consequence of \lemref{UsefulLemma}(iii). 
For, if $a^k\notin P_F$ for some $k\in\N$ then there is some $x\in F$ 
for which $x\leq a^k$, and hence $x=a^m$ for some $m\leq k$. (Here $m>0$ 
since the assumption $F\notin\Cal L(\Omega_\Cal E)$
implies that $e\notin F$). Since $F$ is finite, this can happen for only a finite
number of $a\in S$, leaving at least one for which $a^k\in P_F$ for all $k\in\N$.
Note that, for $a,b\in S$ and $m,n>0$, $a^m\neq b^n$ unless $a=b$, 
again by \lemref{UsefulLemma}(iii). 

In the case that $S$ is finite we have a finite set of lower bounds for $P$.
By \lemref{largest-invariant-set}, the fact that $F\notin\Cal L(\Omega_\Cal E)$ means
that there is some directed hereditary set $\omega$ with no maximal element 
that fails to satisfy the relation $F$.
By invariance, we may choose a translate of $\omega$ and suppose that $\omega\cap F=\emptyset$.
Since $\omega$ is hereditary, $\omega\cap P$ must be contained in $P_F$, 
and since it has no maximal element it must contain an increasing sequence 
$(y_i)_{i\in\N}$ as required. 
\medskip

\paragraph{\bf Step 2.} \emph{The set $P_F$ contains at least 
one element $y$ such that $C(y)$ is a proper subset of $S$.}
\smallskip

It follows from \lemref{UsefulLemma}(ii) that $\ell(u\bsl v)\leq\ell(v)$ for any 
$u,v\in P$. Moreover, if $x,y,z\in P$ we have
\[
z\bsl(y\bsl x)  = z\inv(z\vee y\inv(y\vee x)) = (yz)\inv(yz\vee y\vee x) = (yz)\bsl x\,,
\]
from which it follows that $\ell(y'\bsl x)\leq\ell(y\bsl x)$ whenever $y\leq y'\in P$,
to see this, simply put $y'=yz$ and apply the inequality $\ell(z\bsl(y\bsl x))\leq \ell(y\bsl x)$. 
Therefore, if $(y_i)_{i\in\N}$ is a strictly increasing sequence, as in Step 1,
then $(\ell(y_i\bsl x))_{i\in\N}$ is a decreasing sequence bounded above 
by $\ell(y_1\bsl x)$, for each $x\in F$. Using this and the fact that $F$ is finite, 
we may find  within the sequence $(y_i)_{i\in\N}$ distinct
elements $y<y'$ such that $\ell(y\bsl x)=\ell(y'\bsl x)$ for all $x\in F$. 
Now let $a$ denote any atom dividing $y^{-1}y'$. Then $a$ does not lie in $C(y)$. 
For, if $a\leq\, y\bsl x$, for some $x\in F$, then $(ya)\bsl x=a\inv(y\bsl x)$ and
so $\ell(ya\bsl x) < \ell(y\bsl x)$, strictly. Since $ya\leq y'$, we would then have 
$\ell(y'\bsl x)\leq\ell(ya\bsl x) < \ell(y\bsl x)$, a contradiction. 
It follows that $C(y)$ is a proper subset of $S$, as required.
\medskip

\paragraph{\bf Step 3.} \emph{There exists $y_0\in P_F$ such that $C(y_0)=\emptyset$.}
\smallskip

Consider, amongst all elements $y\in P_F$ for which $C(y)\neq S$, 
an element $y_0$ which minimises the function $L(y):=\Sigma_{x\in F} \ell(y\bsl x)$.
We claim that $C(y_0)=\emptyset$  as required (equivalently $L(y_0)=0$). 
Suppose otherwise.
Then, since $(G,P)$ is graph-irreducible and $C(y_0)$ is a proper nonempty 
subset of $S$, there exist $a\in S\setminus C(y_0)$ and $c\in C(y_0)$ such 
that $a\vee c=\infty$. Set $y=y_0a$, and observe that 
$y\in P_F$. If not, we would have $x\leq y_0a$, for some $x\in F$, and so
$y_0\bsl x = a$  contradicting $a\notin C(y_0)$. 
Now, by the argument used in Step 2 (based on Lemma \ref{UsefulLemma}(ii)), 
we have $\ell(y\bsl x)\leq \ell(y_0\bsl x)$ for each $x\in F$, and so $L(y)\leq L(y_0)$. 
However, since $c\in C(y_0)$, there exists $x_c\in F$ such that $c\leq\, y_0\bsl x_c$.
Because $a\vee c=\infty$, it follows that $y\vee x_c=\infty$. But then 
$\ell(y\bsl x_c)=0$, while $\ell(y_0\bsl x_c)\geq 1$. So, in fact, $L(y)<L(y_0)$. 
Since this contradicts the choice of $y_0$, we must have $C(y_0)=\emptyset$ 
and therefore $F\notin\Cal F$. It follows that  $\partial \Omega = \Omega_\Cal E$, 
which completes the proof of the theorem.
\end{proof}

\begin{example} \label{gammaopexample}
In the case of a right-angled Artin group the different notions of irreducibility coincide.
Given a simplicial graph $\Gamma$, we denote by $\Gamma^{\rm opp}$ its opposite,
or complementary, graph. This is the graph with the same vertex set $S$ but 
where $\{a,b\}$ is an edge of $\Gamma^{\rm opp}$ if and only it is {\em not} 
an edge of $\Gamma$. It is clear that $\Ga^{\rm opp}$ is exactly the $\vee$-graph for
the right-angled Artin group $(A_\Gamma,A_\Gamma^+)$, since the atom set for $A_\Gamma^+$
is precisely the set $S$ of standard generators. Thus the right-angled Artin
group $(A_\Gamma,A_\Gamma^+)$ is graph-irreducible if and only if $\Gamma^{\rm opp}$
is connected. Moreover, any right-angled Artin group decomposes canonically into a direct 
sum of graph-irreducible factors corresponding to the connected 
components of $\Gamma^{\rm opp}$. 
Any direct factor $(\Z,\N)$ of $(A_\Gamma,A_\Gamma^+)$ corresponds to an 
isolated point of $\Gamma^{\rm opp}$, so $A$ has trivial centre if and only if 
$\Gamma^{\rm opp}$ does not have isolated points. 
\end{example}

\begin{corollary}\label{IrredArtin}
If $(A,A^+)$ is a (graph-)irreducible right-angled Artin 
group with standard generating set $S$, then
\[
 \partial\Omega =\Omega_\Cal E =\left\{
\begin{aligned}
\Omega_{\{ S\}} \hsthree &\text{ if $S$ is finite, }\\
\Omega\ \  \hsfour &\text{ if $S$ is infinite. }
\end{aligned}\right.
\] 
\end{corollary}

\begin{proof}
The abelianisation homomorphism $\ell:(A,A^+)\to(\Z,\N)^{|S|}$ is a controlled map, by
\corref{coramenable}. We may therefore apply \proref{IrredCase} to conclude 
that $\partial\Omega =\Omega_\Cal E$. The rest follows from Remark~\ref{partialomegaremark}
and the fact that $S$ is the set of atoms of $A^+$ and therefore the minimal set of 
lower bounds for $A^+$, by \lemref{UsefulLemma}(i). 
\end{proof}

When the Artin group in \corref{IrredArtin} has infinitely many generators,
the boundary quotient  $C(\partial \Omega) \rtimes A$ 
coincides with the universal Toeplitz C*-algebra $C^*(A,A^+)$, which has an
elegant presentation in terms of generators and relations, see \cite[Theorem 24]{CL}. 
It is also possible to give a similar presentation of the boundary quotient of $C^*(A,A^+)$
in the more general case when $A$ is assumed to have trivial centre. 
We underline that the resulting theorem about
the simplicity and pure infiniteness of the C*-algebra 
with the given presentation can be stated with no reference to Artin groups.

\begin{theorem}\label{presentationboundaryquotient}
Suppose $\Gamma$ is a simplicial graph with set of vertices $S$ (finite or infinite)
such that $\Gamma^{\rm opp}$ has no isolated vertices. Then the
 universal C*-algebra with generators $\{V_s: s\in S\}$ subject to the relations 
 \begin{enumerate}
\item $V_s^* V_s = 1 \ \quad \text{ for each } s\in S$; 
\medskip\item
$V_s V_t = V_t V_s \ \  \text{ and } \ \  V_s^* V_t = V_t V_s^*$ \quad 
if $s$ and $t$ are adjacent in  $\Gamma$;
\medskip\item 
$V_s^* V_t = 0\ \quad $ if $s$ and $t$ are distinct and not adjacent in $ \Gamma$;
\medskip\item 
$
\prod_{s\in S_\lambda} (I - V_s V_s^*) = 0$ 
for each $S_\lambda \subset S$ spanning a
finite connected component of $\Gamma^{\rm opp}$,
\medskip
\end{enumerate}
 is canonically isomorphic to the boundary quotient for $(A_\Gamma, A_\Gamma^+)$
  and is purely infinite and simple.
\end{theorem}

\begin{proof}
Since the boundary  quotient of the right-angled Artin group $(A,A^+)=(A_\Gamma, A^+_\Gamma)$ is purely 
infinite and simple by \corref{TopFreeArtin},  it suffices to prove the first assertion.

The first three sets of relations are a presentation of $C^*(A,A^+)=C(\Omega)\rtimes A$, 
so there is a canonical C*-algebra homomorphism
 that sends the standard generator $i(s)$ of $C^*(A,A^+)$ to $V_s$, for each $s\in S$.
In order to conclude that this map induces an isomorphism at the level of the boundary quotient, 
it suffices to show that an element $\omega\in\Omega$ lies in the boundary spectrum
$\partial\Omega$ if and only if each of
the relations given in (4) is satisfied at $\omega$.

The Artin group $(A, A^+)$ has a canonical decomposition
as a direct sum of graph-irreducible Artin groups 
\[
(A,A^+)=\bigoplus_{\lambda\in\Lambda}(A_\lambda, A^+_\lambda)
\]
corresponding to the decomposition of  $\Gamma^{\rm opp}$ into its connected components. 
Here the indexing set $\Lambda$ may be finite or countably infinite. We let $S_\lambda$ denote the 
standard generating set for the direct factor $(A_\lambda, A^+_\lambda)$. Thus, viewed as a subset of
$S$, each $S_\lambda$ spans a connected component of $\Gamma^{\rm opp}$.

Writing $\Omega$ for the Nica spectrum of $(A,A^+)$ and $\Omega_\lambda$ for that
 of $(A_\lambda, A^+_\lambda)$, for each $\lambda\in\Lambda$, we observe that
if $\omega\in\Omega$ then $\omega=\bigoplus\limits_{\lambda\in\Lambda}\omega_\lambda$ 
where $\omega_\lambda\in\Omega_\lambda$ for each $\lambda\in\Lambda$. Also, $\omega$ 
is clearly a \emph{maximal} element of $\Omega$ precisely when each $\omega_\lambda$ 
is maximal as an element of $\Omega_\lambda$. Moreover a sequence $(\omega_{i})_{i\in\N}$
converges to $\omega$ in $\Omega$ if and only if the sequences 
$((\omega_{i})_\lambda)_{i\in\N}$ converge to $\omega_\lambda$ in 
$\Omega_\lambda$ for all $\lambda$.
Therefore, by \lemref{minimalityofboundary}, we have
 $\omega\in\partial\Omega$ if and only if $\omega_\lambda\in\partial\Omega_\lambda$
for every $\lambda\in\Lambda$.
However, by \corref{IrredArtin}, we have $\omega_\lambda\in\partial\Omega_\lambda$
 if and only if either $S_\lambda$ is infinite, or $\omega_\lambda$ satisfies the
elementary relation $S_\lambda$ (as an element of $\Omega_\lambda$). 
Finally we observe that $\omega\in\Omega$ satisfies the relation $S_\lambda$ if and only
if $\omega_\lambda$ does. It follows that an element $\omega\in\Omega$ lies
in $\partial\Omega$ if and only if it satisfies each of the relations given 
in (4).

\end{proof}

\begin{remark}
When  no connected component of $\Gamma^{\rm opp}$ is finite, the relations 
in \thmref{presentationboundaryquotient} (4) never arise and $C^*(A,A^+)$ itself is purely 
infinite and simple. Note that since such an Artin group satisfies 
condition (ii) of \cite[Lemma 5.2]{purelinf}, 
this assertion already follows from \cite[Theorem 5.4]{purelinf}.
  When $S$ is finite but not a singleton, and $\Gamma^{\rm opp}$ is connected,
the set of relations (4) reduces to a single relation:
\begin{equation}
\prod_{s\in S} (I - V_s V_s^*) = 0;
\end{equation}
At the level of $C^*(A,A^+)$, the ideal generated by the projection 
$\prod_{s\in S} (I - V_s V_s^*) $ is isomorphic to the compact 
operators on a separable Hilbert space by \cite[Proposition 6.3]{nica}.
\end{remark}


\section{Direct sums of quasi-lattice orders and the Nica spectrum}\label{SectProducts}


Recall that any right-angled Artin group $(A,A^+)$ is canonically a direct
sum of graph-irreducible right-angled Artin groups.
We wish to reduce the question of topological freeness of $G$ on 
closed invariant subsets of $\Omega$ to the analogous question for each 
irreducible factor of the quasi-lattice order $(G,P)$.
We first introduce some basic notions and notations 
regarding direct products of partial actions.

\begin{notation}
In general, we shall write $(G,X)$ to denote a partial action of a group $G$
by partial homeomorphisms of a locally compact space $X$. Given a family 
$\{(G_\lam,X_\lam) : \lam\in\Lam\}$ of such partial actions, we define the
\emph{direct product of partial actions}
\[
(G,X)= \prod_{\lam\in\Lam}(G_\lam,X_\lam)
\]
where $G=\oplus G_\lam$, $X=\prod X_\lam$ and where, for $g=(g_\lam)_{\lam\in\Lam}\in G$, 
we define the range of $g$ on $X$ by $U_g=\prod U_{g_\lam}$
and define the action of $g$ on a point $x=(x_\lam)_{\lam\in\Lam}\in U_{g\inv}$ by 
$g(x)=(g_\lam(x_\lam))_{\lam\in\Lam}$.
\end{notation}

\begin{lemma}\label{TopFreeProducts}
Suppose that $(G,X)= \prod_{\lam\in\Lam}(G_\lam,X_\lam)$ is a direct product of 
partial actions. Then $(G,X)$ is topologically free if and only if
$(G_\lam,X_\lam)$ is topologically free for every $\lam\in\Lam$. 
\end{lemma}

\begin{proof}
We first note that the statement that $(G,X)$ is topologically free is equivalent to
saying that if $g\in G$ fixes pointwise a nonempty open set $U\subset X$ then
$g=e$.

Suppose firstly that $(G_\lam,X_\lam)$ is topologically free for every $\lam\in\Lam$. 
If $g=(g_\lam)_{\lam\in\Lam}\in G$ then $g_\lam=e$ for all $\lam$ but a finite set 
$\{1,..,k\}\subset\Lam$. Writing 
\[
X=X_1\times\dots\times X_k\times X'\,, \hsfour \text{ where }
X'=\prod_{\lam\notin\{1,..,k\}} X_\lam\,,
\]
we have $g=(g_1,..,g_k,e)$. If $U\subset X$ is a nonempty open set then it
contains a subset of the form $U_1\times\dots U_k\times U'$ where each $U_i$ is a nonempty
open subset of $X_i$ (and $U'$ a nonempty open subset of $X'$). Now if $g$ fixes $U$
pointwise we have that $g_i$ fixes $U_i$ pointwise, for each $i=1, \ldots ,k$. But then since
each $(G_i,X_i)$ is topologically free we deduce that each $g_i=e$ and hence $g=e$.
This shows that $(G,X)$ is topologically free. 

On the other hand, if we suppose that $(G,X)$ is topologically free then
each factor $(G_\mu,X_\mu)$ must be also. For if $g_\mu$ fixes pointwise a nonempty 
open set $U_\mu\subset X_\mu$ then setting $g_\lam=e$ for all $\lam\neq\mu$, the group element 
$g=(g_\lam)_{\lam\in\Lam}$ fixes pointwise the open set $p_\mu\inv(U_\mu)\subset X$, 
where $p_\mu:X\to X_\mu$ denotes the canonical projection.
Since $(G,X)$ is topologically free, $g = e$ and hence $g_\mu$ is trivial.
\end{proof}

\begin{lemma}\label{TopFreeUnions}
Given a partial action $(G,\Omega)$ we suppose that $\{\, X_i\,:\, i\in I\,\}$
is a collection of closed $G$-invariant subsets of $\Omega$.
If the restricted action $(G,X_i)$ is topologically free for every $i\in I$ 
then the action $(G,\bar Z)$ is topologically free, where $\bar Z$ denotes
 the closure of  the set $Z=\bigcup\limits_{i\in I} X_i$.
\end{lemma}

\begin{proof}
If $Y\subset\Omega$ is a closed invariant subset then $(G,Y)$ is topologically free
if and only if, for every nontrivial $g\in G$ and open set $U\subset\Omega$,
we have that if $g$ fixes $U\cap Y$ then
$U\cap Y=\emptyset$. Suppose that each $(G,X_i)$ is topologically free.
Then, since 
$U\cap Z=U\cap(\bigcup\limits_{i\in I} X_i)=\bigcup\limits_{i\in I}(U\cap X_i)$, 
we have
\[
\begin{aligned}
g\text{ fixes }U\cap \bar Z 
&\implies g\text{ fixes }U\cap Z
\implies g\text{ fixes }U\cap X_i\text{ for all }i\in I \\
&\implies U\cap X_i=\emptyset\text{ for all }i\in I 
\implies U\cap Z=\emptyset\implies U\cap \bar Z=\emptyset\,,
\end{aligned}
\]
for every nontrivial $g\in G$ and open set $U\subset \Omega$. 
\end{proof}

Suppose now that $(G,P)=\oplus_{\lam\in\Lam} (G_\lam,P_\lam)$ is a direct sum
of quasi-lattice orders $(G_\lam,P_\lam)$. As usual we write $\Omega$ for the Nica spectrum
associated to $(G,P)$ and $\Omega_\lam$ for the Nica spectrum of the summand
$(G_\lam,P_\lam)$, for each $\lam\in\Lam$. We then have the following:

\begin{lemma}\label{NicaProd}
The Nica spectrum $\Omega$ is canonically homeomorphic to 
$\prod_{\lam\in\Lam} \Omega_\lam$, where the canonical partial action of $G$ on $\Omega$
is just the product of the canonical partial actions of the direct summands: namely,
for $g=(g_\lam)_{\lam\in\Lam}$ and $\omega=(\omega_\lam)_{\lam\in\Lam}$ we have 
$g(\omega)=(g_\lam(\omega_\lam))_{\lam\in\Lam}$.
\end{lemma}

\begin{proof}

Let $p_\lam:G\to G_\lam$ and $i_\lam:G_\lam\to G$ denote the canonical 
projection and inclusion for each summand $G_\lam$ of $G$. The projection 
$p_\lambda$ naturally induces a continuous map (also denoted $p_\lambda$) 
from $\Omega$ onto $\Omega_\lam$, for each $\lam\in\Lam$. 
To check continuity, observe that if $t\in G_\lam$ and $H$ is a finite 
subset of $P_\lam$ determining a basic open set $V(t,tH)\subset \Omega_\lam$
then $p_\lam\inv(V(t,tH))= V(i_\lam(t),i_\lam(tH))$. 
(Note that if $\omega\in\Omega$ and $x\in G_\lam$ then  
$i_\lam(x)\in\omega$ if and only if $x\in p_\lam(\omega)$).
It follows that the product map 
\[
f:=(p_\lam)_{\lam\in\Lam}:\Omega\to\prod_{\lam\in\Lam} \Omega_\lam\,,
\]
is continuous.
 It is clear that $f$ is equivariant 
with respect to the canonical partial actions, as claimed in the Lemma.
Next we shall show that $f$ is a homeomorphism.

Given ${\underline\omega}=(\omega_\lam)_{\lam\in\Lam}\in\prod_{\lam\in\Lam}\Omega_\lam$ 
we define $\what{\underline\omega}$ to be the directed hereditary closure of the union
of the subsets $i_\lam(\omega_\lam)$ in $G$. We claim that 
\[
\what{\underline\omega}=\{ g\in G\,:\, p_\lam(g)\in\omega_\lam 
\text{ for all } \lam\in\Lam\}\,.
\]
It is easy to check that the right hand side is hereditary, directed and contains 
all the $i_\lambda(\omega_\lambda)$, and so it contains $\what{\underline\omega}$.
Conversely, suppose $\beta$ is a directed hereditary subset of $G$ containing every 
$i_\lambda(\omega_\lambda)$. If $g\in G$ satisfies $p_\lambda(g) \in \omega_\lambda$,
then $ i_\lambda(p_\lambda(g)) \in i_\lambda(\omega_\lambda) \subset \be$.
Since $g$ is the least common upper bound of its coordinates
$ i_\lambda(p_\lambda(g))$ and $\be$ is directed, then 
 $g \in \be$, finishing the proof of the claim.

It follows that, for any ${\underline\omega}$ as above, 
$p_\lam(\what{\underline \omega})=\omega_\lam$ for
each $\lam\in\Lam$, and so $\underline \omega = f(\what{\underline\omega})$.
Hence $f$ is surjective.
 On the other hand, given $\omega\in\Omega$, one easily checks that 
$\what{(p_\lam(\omega))_{\lam\in\Lam}}=\omega$, 
from which it follows that $f$ is injective.
This shows that $f$ is a continuous bijection between compact 
Hausdorff spaces,
hence, by a standard argument, its inverse is also continuous, finishing the proof.
\end{proof}

\begin{proposition}\label{TopFreeReduction}
Suppose that $(G,P)=\oplus_{\lam\in\Lam} (G_\lam,P_\lam)$ is a direct sum of
quasi-lattice orders $(G_\lam,P_\lam)$, and let $\Omega$ and $\Omega_\lam$, $\lam\in\Lam$,
denote the associated Nica spectra. Then the following are equivalent:
\begin{itemize}
\item[(i)] The canonical partial action of $G$ on $\Omega$ is topologically free
on every closed invariant subset.
\item[(ii)] For each $\lam\in\Lam$, the canonical partial action 
of $G_\lam$ on $\Omega_\lam$ is topologically free on every closed invariant subset.
\end{itemize}
\end{proposition}

\begin{proof}
Recall that by Lemma \ref{NicaProd} the canonical action $(G,\Omega)$ may be written as
a product of actions 
\[
(G,\Omega)=\prod_{\lam\in\Lam} (G_\lam,\Omega_\lam)\,.
\]
Let $\omega\in\Omega$ and let  $X(\omega)$ be the 
smallest closed invariant subset
 of $\Omega$ containing $\omega$, i.e{.}{} the intersection 
of all closed invariant subsets containing $\omega$.
If $\omega=(\omega_\lam)_{\lam\in\Lam}$, then we may write $X_\lam(\omega_\lam)$
for the smallest closed invariant subset of $\Omega_\lam$ containing $\omega_\lam$.
One easily observes that $X(\omega)=\prod_{\lam\in\Lam}X_\lam(\omega_\lam)$.

It follows from the above observation that any closed 
invariant set $X\subset \Omega$ may be expressed as 
the closure of  a union of products of
closed invariant sets belonging to the factors $\Omega_\lam$,
\[
X=\overline{\bigcup_{\omega\in X}X(\omega)}\ =\ 
\overline{\bigcup_{\omega\in X}\ \prod_{\lam\in\Lam}X_\lam(\omega_\lam)}\,\,.
\]
The statement that  (ii) implies (i) now follows by Lemmas \ref{TopFreeProducts} 
and \ref{TopFreeUnions}. The converse, that (i) implies (ii), may be easily deduced
from Lemma \ref{TopFreeProducts} and the fact that if $X_0$ is a closed invariant subset
of $\Omega_{\lam_0}$ then $X_0\times\prod_{\lam\neq\lam_0}\Omega_\lam$ is a closed 
invariant subset of $\Omega$. 
\end{proof}

A closed invariant subset of $\Omega$ shall be called a \emph{component} subset 
if it is of the form $X=\prod_{\lam\in\Lam} X_\lam$ where $X_\lam =\Omega_\lam$ for all
but at most one factor. Equivalently, a component subset $X$ of $\Omega$
 is the spectrum of a family of elementary 
relations coming from a single summand $(G_{\lam_0},P_{\lam_0})$ of $(G,P)$. 
The proof of Proposition \ref{TopFreeReduction} above rests on the observation
that every closed invariant set of $\Omega$ can be expressed as 
 the closure of a union of intersections of component subsets. 

We would like to interpret the above decomposition of closed invariant sets in terms
of defining relations. We first observe that, quite generally, 
for $H,K\in\Cal Q$ and $\Cal R\subset\Cal Q$,
one has
\[ 
\Omega_{\{H\cup K\}}=\Omega_{\{ H\}}\cup\Omega_{\{ K\}}\qquad \text{ and }\qquad 
\Omega_{\Cal R}= \bigcap_{R\in\Cal R}\Omega_{\{ R\}}.
\]
We shall say that a closed invariant subset $X$ of $\Omega$ is \emph{principal} 
if it can be written as the spectrum of a single elementary relation, 
namely $X=\Omega_{\{H\}}$ for some $H\in\Cal Q$. A \emph{principal component subset}
is therefore a closed invariant set of the form $\Omega_{\{ H\}}$ where the elementary
relation $H$ lies in a single irreducible summand of $(G,P)$.

The equality $\Omega_{\Cal R}= \bigcap_{R\in\Cal R}\Omega_{\{ R\}}$
stated above, together with Lemma \ref{ClosedSetsRelations},
implies that every closed invariant subset of $\Omega$ is expressible 
as the intersection of a countable collection of principal closed invariant sets.

In particular, \emph{any component subset may be written as the intersection of a 
countable collection of principal component subsets}. 
We now have:

\begin{proposition}\label{PrincipalComponents}
Suppose that $(G,P)=\oplus_{\lam\in\Lam} (G_\lam,P_\lam)$ is a direct sum of
quasi-lattice orders $(G_\lam,P_\lam)$, and let $\Omega$ and $\Omega_\lam$, $\lam\in\Lam$,
denote the associated Nica spectra. Suppose moreover that, for each $\lam\in\Lam$, the 
component Nica spectrum $\Omega_\lam$ has only finitely many closed invariant subsets. 
Then, every closed invariant subset of $\Omega$ may be written as a countable intersection
of finite unions of principal component subsets. In other words, every closed invariant
subset of $\Omega$ is the spectrum $\Omega_{\Cal R}$ of a set $\Cal R$ of elementary relations 
such that each relation $R\in\Cal R$ is the union of finitely many subsets, each 
belonging to a single direct summand $P_\lam$ of $P$. 
\end{proposition}
   
\begin{proof}
By the preceding remarks it clearly suffices to prove the statement for any principal 
closed invariant set $X=\Omega_{\{ H\}}$, for $H\in\Cal Q$. 
Since $H$ is a finite set of elements
of $P=\oplus_{\lam\in\Lam} P_\lam$, there is a finite collection 
$\{\lam_1,..,\lam_m\}\subset\Lam$ such that $H$ is contained in 
the factor $P'=P_{\lam_1}\times \dots\times P_{\lam_m}$. We may write
\[
\Omega =\Omega'\times\Omega''\hsthree \text{ where }\hsthree 
\Omega'=\Omega_{\lam_1}\times\dots\times\Omega_{\lam_m} \hsthree \text{ and }\hsthree
\Omega''=\prod_{\lam\notin\{\lam_1,..,\lam_m\}}\Omega_\lam\,.
\]
Then $X=\Omega_{\{ H\}}$ has the form $X'\times\Omega''$, 
where $X'$ is a closed invariant subset of $\Omega'$. By our previous remarks, $X'$
may be expressed as  the closure of  a union of intersections of 
principal component subsets of $\Omega'$.
However, since we suppose that every component $\Omega_\lam$ has only finitely many
closed invariant subsets, it follows that $\Omega'$ has only finitely many component
subsets and, in particular, only finitely many principal component subsets.
Thus any closed invariant subset of $\Omega'$, and in particular $X'$, can in fact 
be written as a \emph{finite} union of \emph{finite} intersections of principal
component subsets. 
Alternatively, by de Morgan's Law, $X'$ may be written as a finite intersection 
of finite unions of principal component subsets. Restoring the 
factor $\Omega''$ everywhere we now have $X$ as a finite intersection of finite 
unions of principal component subsets. 
\end{proof}


\section{The lattice of ideals of $C^*(A,A^+)$ 
for a right-angled Artin group}\label{SectIdealStruct}


In \cite{CL} it was proved that, in the case of a right-angled Artin group $(A,A^+)$,
the Toeplitz algebra $\Cal T(A,A^+)$ is universal for isometric representations 
satisfying Nica's covariance condition. That is, $\Cal T(A,A^+)\cong C^*(A,A^+)$. 
By combining results of the preceding sections we are now able to completely 
describe the structure of the lattice of ideals of this C*-algebra 
when $(A,A^+)$ has no direct $(\Z,\N)$ factor (equivalently, when $A$ has 
trivial centre, or $\Gamma^{\rm opp}$ has no isolated points).

\begin{proposition}\label{IdealsToSets}
Let $(A,A^+)$ denote a right-angled Artin group with the standard quasi-lattice order
and let $\Omega$ denote the associated Nica spectrum.
Suppose that $(A,A^+)$ has no direct $(\Z,\N)$ factor.
Then the canonical partial action of $A$ on $\Omega$ is 
amenable and topologically free on every closed invariant subset, and 
the map $F \mapsto \< C_0(\Omega \setminus F)\> $ is an inclusion-reversing isomorphism between the 
lattice of closed invariant subsets of $\Omega$ and the lattice of 
ideals in $C^*(A,A^+)$.
\end{proposition}

\begin{proof} 
By \corref{coramenable}, the restriction of 
the canonical partial action of $A$ to each closed invariant subset of $\Omega$ is amenable.
Let $\Gamma$ be the graph defining $A$, and recall from Example \ref{gammaopexample} 
that $(A,A^+)$ decomposes canonically as a direct sum of graph-irreducible 
right-angled Artin groups, corresponding to the connected components of the opposite 
graph $\Gamma^{\rm opp}$. 
In the irreducible case, by \corref{IrredArtin}, there are at most
two closed invariant sets, namely $\Omega$ and $\partial\Omega$.
By \cite[Proposition 6.6]{ELQ} the action is 
always topologically free on $\Omega$, and topological freeness on $\partial\Omega$
is proved in \corref{TopFreeArtin} under the assumption that $\Gamma $ has no $(\Z,\N)$ factor. 
It follows, by Proposition \ref{TopFreeReduction}, that the partial action is 
topologically free on all closed invariant sets in the general case, and an application of 
\cite[Theorem 3.5]{ELQ}, as restated in \thmref{IdealStructure}, finishes the proof.
\end{proof}

Let $(A,A^+)$ be a right-angled Artin group with Nica spectrum $\Omega$.
If $\Cal R$ is a set of elementary relations in $A^+$ 
or, more precisely, if $\Cal R$ is a collection of finite subsets  of $A^+$,
then we may speak of the \emph{ideal generated by $\Cal R$} to mean 
the ideal of $C^*(A,A^+)$ generated by the elementary functions 
$f_H :=\prod_{h\in H} (1 - 1_h)$ for $H\in \Cal R$. This ideal 
corresponds  to the closed invariant subset $\Omega_\Cal R$
under the correspondence of Proposition \ref{IdealsToSets}. We shall say that an 
elementary relation $H$ \emph{is a consequence of} the elementary relation $K$
if $f_H =0$ whenever $f_K =0$, that is, 
if $\Omega_{\{ H\}}\supseteq\Omega_{\{ K\}}$. 

\begin{definition}
Let $(A,A^+)=\oplus_{\lam\in\Lam} (A_\lam,A^+_\lam)$ denote the canonical 
decomposition of a right-angled Artin group into its irreducible factors.
For each $\lam\in\Lam$ we shall write $S_\lam$ for the standard generating set
of the summand $(A_\lam,A^+_\lam)$, 
and write $\Lam_f=\{ \lam\in\Lam : S_\lam \text{ finite}\,\}$ for the
set indices associated to the finitely generated irreducible direct factors of $(A,A^+)$,
corresponding to the finite connected components of $\Gamma^{\rm opp}$. 
By a \emph{basic} relation we mean any elementary relation of the form 
\[
S_B:=\cup_{\lam\in B} S_\lam\qquad \text{for } B \text{ a finite subset of } \Lam_f\,. 
\]
We denote by $\Cal P(\Lam_f)$ the Boolean lattice of finite subsets of $\Lam_f$
ordered by \emph{reverse} inclusion. 
Note that $\Cal P(\Lam_f)$ has a maximal element, namely the empty set $\emptyset$,
while the basic relation $S_\emptyset =\emptyset$, interpreted as the empty product
$f_\emptyset =1\in C_0(\Omega)$,  generates the whole algebra $C^*(A,A^+)$ as an ideal. 

In the usual combinatorial terminology, a hereditary subset 
of a partially ordered set $\Cal P$ is referred to as an \emph{ideal} of $\Cal P$. 
We include the empty set as an ideal, the empty ideal. 
The ideals of a partially ordered set form a lattice, where 
meets and joins are given by taking intersections and unions of ideals respectively.
An ideal of $\Cal P$ is said to be generated by a subset of $\Cal P$ if it is the 
smallest ideal containing the given subset. 
 Finally, we observe that the ideals of 
$\Cal P(\Lam_f)$ form a \emph{complete} lattice, with maximal ideal the ideal generated by 
$\emptyset$, and minimal ideal the empty one. 
\end{definition}

\begin{theorem}\label{LatticeOfIdeals}
Let $(A,A^+)=\oplus_{\lam\in\Lam} (A_\lam,A^+_\lam)$ denote the canonical 
decomposition of a right-angled Artin group into its irreducible factors, and 
suppose that no factor of this decomposition is isomorphic to $(\Z,\N)$.
For each ideal $\Cal B \in \Cal P(\Lam_f)$ let $\phi(\Cal I) $ denote the ideal 
$ \langle  \prod_{s \in S_B} (1 - 1_s) : B \in \Cal B\rangle$, generated by the basic relations
 $S_B $ for $ B \in \Cal B$.
Then $\phi$ is an order isomorphism of the lattice of ideals of the Boolean lattice $\Cal P(\Lam_f)$ to
the lattice of ideals of the C*-algebra $C^*(A,A^+)$, which maps
the principal ideal of $\Cal P(\Lam_f)$ generated by the finite set $B \subset \Lambda_f$
to  the principal ideal generated by the basic relation $S_B$.
\end{theorem}

\begin{proof}
We first observe that if  $\Cal B \subset  \Cal C$  in $ P(\Lam_f)$,  then clearly 
$\phi(\Cal B) \subset \phi(\Cal C)$, so $\phi$ is an order preserving map from 
the lattice of ideals of $\Cal P(\Lam_f)$ into the lattice of ideals of $C^*(A,A^+)$.
Also notice that if $\Cal B = \langle B \rangle$, then $\phi(\Cal B) = \langle S_B \rangle$.

Let $\Omega^{(\lam)}$ denote the Nica spectrum for the 
irreducible factor $(A_\lam,A^+_\lam)$. By Corollary \ref{IrredArtin} we have that 
$\partial\Omega^{(\lam)}=\Omega^{(\lam)}_{\{S_\lam\}}$ if $\lam\in\Lam_f$ and 
$\partial\Omega^{(\lam)}=\Omega^{(\lam)}$ otherwise.
These being the only closed invariant sets in each factor,
any finite union of principal component subsets is necessarily of the form
$\Omega_{\{ S_B\}}$ for some basic relation $S_B$. (Recall the equation 
$\Omega_{\{H\cup K\}}=\Omega_{\{ H\}}\cup\Omega_{\{ K\}}$). 
It follows from Proposition \ref{PrincipalComponents} that 
 every closed invariant subset of $\Omega$ 
is of the form $\Omega_\Cal B$ for $\Cal B$ a collection of basic relations.
In order to see that $\phi$ is surjective, let $\Cal I$ be an ideal in $C^*(A,A^+)$. 
By \thmref{IdealsToSets} $\Cal I$ is generated by a closed invariant subset of $\Omega$, 
which is of the form $\Omega_\Cal B$. It is easy to check that $\phi(\Cal B) = \Cal I$.

It remains to show that 
$\phi(\Cal B)\leq\phi(\Cal C)$ implies $\Cal B \leq \Cal C$, 
for all ideals $\Cal B,\Cal C$ of $\Cal P(\Lam_f)$; clearly this implies injectivity of $\phi$ and 
will complete the proof.
Given  $B\in \Cal P(\Lam_f)$, we may write 
\[
(G_B,P_B)=\prod_{\lam\in B} (G_\lam,P_\lam)\hsthree \text{ and }\hsthree
(G^\perp_B,P^\perp_B)=\prod_{\lam\in\Lam\setminus B} (G_\lam,P_\lam)\,,
\]
and define the directed hereditary set
\[
\omega_B = P_B^{-1}\times G^\perp_B\ \subset\ G=G_B\times G^\perp_B \,.
\]
Let $C$ be a finite subset of $\Lam_f$, and observe that $\omega_B$ satisfies 
the basic relation $S_C$ if and only if $C$ is not a subset of $B$.
In particular, $\omega_B$ does not satisfy $S_B$.  It follows that if the
relation $S_B$ is a consequence of some set of basic relations $\Cal R$ then
there is some $S_C\in\Cal R$ with $C\subseteq B$, for otherwise we would 
have $\omega_B\in\Omega_{\Cal R}$ but $\omega_B\notin\Omega_{\Cal R\cup\{ f_B\}}$. 

Suppose  now that the ideal $\Cal B$  
is generated by the set $\{ B_i : i\in I\}$ of elements
of $\Cal P(\Lam_f)$, and similarly, that the ideal $\Cal C$ 
is generated by the set $\{ C_j : j\in J\}$.
Then $\phi(\Cal B)$ is generated by the collection of basic relations 
$\{S_{B_i} : i\in I\}$ and $\phi(\Cal C)$  by the collection  
$\{S_{C_j} : j\in J\}$. Suppose that $\phi(\Cal B)\leq\phi(\Cal C)$.
Then, for each $i\in I$, the basic relation $S_{B_i}$ is a consequence 
of the relations $\Cal C$ and so, by the previous argument, there is 
some $j\in J$ such that $C_j\subseteq B_i$, or rather $B_i\leq C_j$. 
It follows that $\Cal I\leq\Cal J$.
\end{proof}

\begin{remark}
Note that each ideal of $\Cal P(\Lam_f)$ is uniquely generated by a set of mutually
incomparable elements (two finite sets $B,C\in\Cal P(\Lam_f)$ are said to be 
incomparable if neither
$B\subseteq C$ nor $C\subseteq B$). This leads to an effective method for performing
computations in the lattice of ideals of  $\Cal P(\Lam_f)$, and hence in the lattice of
ideals of $C^*(A,A^+)$ for any right-angled Artin group with trivial centre.
\end{remark}

Finally, from this point of view it is clear that, for a right-angled Artin group $A$ with trivial
centre, the maximal proper ideal of $C^*(A,A^+)$ is generated by the basic relations
$\{ S_B\ :\ B\text{ is a singleton }\}$, that is the elementary relations 
$\{ S_\lam\ :\ \lam\in\Lam_f \}$ corresponding to the finite connected components of 
$\Gamma^{\rm opp}$. This leads immediately to the set of relations given in 
\thmref{presentationboundaryquotient} (4) for the presentation of the boundary quotient. 
More generally, we have the following statement: 

\begin{corollary}\label{quotientpresentation}
Suppose $\Gamma$ is a simplicial graph with set of vertices $S$ (finite or infinite)
such that $\Gamma^{\rm opp}$ has no isolated vertices. Let $C^*(\Gamma)$ denote
the universal C*-algebra with generators $\{V_s: s\in S\}$ subject to the relations 
\begin{enumerate}
\item $V_s^* V_s = 1$ \ for each $s\in S$, 
\item
$V_s V_t = V_t V_s \ \  \text{ and } \ \  V_s^* V_t = V_t V_s^*$ \ 
if $s$ and $t$ are adjacent in  $\Gamma$,
\item 
$V_s^* V_t = 0$\  if $s$ and $t$ are distinct and not adjacent in $ \Gamma$.
\end{enumerate}
Then each quotient of $C^*(\Gamma)$ is obtained by imposing 
a further collection of relations of the form
\begin{enumerate}
\item[(R)] $\prod_{s\in S_i} (I - V_s V_s^*) = 0 $ for $i\in I$, where each $S_i \subset S$ 
spans a finite union of finite connected components of $\Gamma^{\rm opp}$. 
\end{enumerate}
Moreover, we may always reduce the above presentation to one in which 
no two sets $S_i$ and $S_j$ used are contained one in the other. 
\end{corollary} 

\begin{example}
(a) Suppose $\Gamma$ has a finite  set of vertices $S$ and no edges, so that $A_\Gamma $ is the free group $\mathbb F_{|S|}$ on $|S|$ generators and $C^*(\Gamma)$ is the Toeplitz-Cuntz algebra $\mathcal T_{|S|}$.
In this case $\Gamma^{\rm opp}$ is connected and the only possible nontrivial collection of relations, 
to be added as in part (R) above, is Cuntz's $\prod_{s\in S}(1-1_s)=  1-\sum_{s\in S} 1_s  = 0$.

(b) Suppose $\Gamma$ consists of two finite 
sets of vertices $A$ and $B$ with no edges from $A$ to $A$ 
or from $B$ to $B$ but with everything in $A$ adjacent to everything in $B$. Then
$A_\Gamma \cong \mathbb F_{ |A|} \times \mathbb F_{|B|}$, and $C^*(\Gamma)$ is the tensor product of the
Toeplitz -Cuntz algebras associated to $A$ and $B$. In this case $\Gamma^{\rm opp} $
has two (finite) connected components $A$ and $B$, and the possible (reduced) sets of additional relations 
are 
\begin{enumerate}
\item[(R1)]  $\prod_{s\in A\cup B}(1-1_s) = 0$
\item[(R2)]  $\prod_{s\in A}(1-1_s)= 0$;
\item[(R3)] $\prod_{s\in B}(1-1_s)= 0$
\item[(R4)] $\prod_{s\in A}(1-1_s) = 0$ and  $\prod_{s\in  B}(1-1_s) = 0$
\end{enumerate}
giving rise to the four proper quotients, with the boundary quotient being given by (R4).

(c) Suppose $\Gamma$ has four vertices $\{a, b, c, d\}$ and three edges $\{a, b\}$ $\{b, c\}$ $\{ c, d\}$ 
(the Artin group of $\Gamma$ is not a direct nor a free product of smaller ones).
Then $\Gamma^{\rm opp}$ is connected, so there is only one nontrivial set of relations
(R): $(1-1_a)(1-1_b)(1-1_c) (1-1_d)$ and there is only one nontrivial quotient of $C^*(\Gamma)$, 
namely the boundary quotient.
\end{example}

Note that the different quotients of the algebra $C^*(\Gamma)$ are distinguished
(as quotients of the same algebra) by the reduced presentations described 
in \corref{quotientpresentation}. However, we do not yet know how to classify these algebras up to
$*$-isomorphism. Nor do we know, for that matter, how to classify the boundary quotients of 
$C^*(\Gamma)$, for different $\Gamma$'s, up to $*$-isomorphism. On the positive side
with respect to this question, 
in ongoing joint work with B. Abadie, we have been able to show in several situations
 (but believe to be true in general) that the order  in $K_0$ of the identity element 
 of the boundary quotient of $C^*(\Gamma)$
 is $|\chi(\Gamma) -1|$, where $\chi (\Gamma)$ is the Euler characteristic of $\Gamma$, viewed as a simplicial complex. 
 
\bigskip
\noindent {\bf Acknowledgment: }
This research was initiated during a visit of J.C. to the SFB 478 at the 
University of M\"unster, and 
continued through visits of M.L. to the IMB at the Universit\'{e} de Bourgogne, Dijon.
We would like to thank those institutions for their support and their hospitality. We 
would also like to thank Beatriz Abadie for several helpful conversations.

\end{document}